\documentclass[12pt,a4paper]{article}

\usepackage[table]{xcolor}
\usepackage{slashed}
\usepackage{xkeyval}
\usepackage{url,amsmath,amssymb,latexsym,pstricks,mathrsfs,comment,amsthm,graphicx,tikz,tikz-cd,enumerate,accents,pgffor,cite,wrapfig,multicol,float}
\usepackage{bibspacing}

\allowdisplaybreaks

\newcommand{\nc}{\newcommand}
\nc{\rnc}{\renewcommand}

\nc{\RP}{\mathcal{RP}}
\nc{\supp}{\operatorname{supp}}

\nc{\OEIS}{}

\usepackage{geometry}  \rnc{\ss}{\smallskip} \nc{\ms}{\medskip}  \nc{\nss}{\vspace{-3mm}} %\nc{\bs}{\bigskip}
\makeatletter

\DeclareMathSymbol{\widehatsym}{\mathord}{largesymbols}{"62}
\newcommand\lowerwidehatsym{%
  \text{\smash{\raisebox{-1.3ex}{%
    $\widehatsym$}}}}
\newcommand\fixwidehat[1]{%
  \mathchoice
    {\accentset{\displaystyle\lowerwidehatsym}{#1}}
    {\accentset{\textstyle\lowerwidehatsym}{#1}}
    {\accentset{\scriptstyle\lowerwidehatsym}{#1}}
    {\accentset{\scriptscriptstyle\lowerwidehatsym}{#1}}
}
\rnc{\widehat}{\fixwidehat}

\begin{document}

%\rnc{\Box}{\includegraphics[width=3.5mm]{like.png}} % FB ``Like'' for end of proof 

\nc{\PnSn}{\P_n\sm\S_n}
\nc{\RPnSn}{\RP_n\sm\S_n}
\nc{\bbR}{\mathbb R}
%\nc{\bbE}{\mathbb E}
\nc{\bbS}{\mathbb S}
\nc{\bbA}{\mathbb A}
\nc{\bbB}{\mathbb B}
\nc{\pre}{\preceq}
\nc{\Sub}{\operatorname{Sub}}
\nc{\MSA}{M(\S_n,\bbA_n)}
\nc{\Th}{\Theta}
\nc{\divider}{\big|}
\nc{\alh}{\widehat{\al}}

%\nc{\uvert}[1]{\fill (#1,2)circle(.2);}
%\rnc{\lvert}[1]{\fill (#1,0)circle(.2);}
\nc{\uvertcol}[2]{\fill[#2] (#1,2)circle(.2);}
\nc{\lvertcol}[2]{\fill[#2] (#1,0)circle(.2);}
\nc{\uvertcols}[2]{\foreach \x in {#1}{ \uvertcol{\x}{#2}}}
\nc{\lvertcols}[2]{\foreach \x in {#1}{ \lvertcol{\x}{#2}}}
\nc{\uverts}[1]{\foreach \x in {#1}{ \uvert{\x}}}
\nc{\lverts}[1]{\foreach \x in {#1}{ \lvert{\x}}}
\nc{\uvws}[1]{\foreach \x in {#1}{ \uvw{\x}}}
\nc{\lvws}[1]{\foreach \x in {#1}{ \lvw{\x}}}
%\nc{\urookdot}[1]{\draw[fill=white] (#1,2)circle(.2);}

\nc{\PBnt}{R^\tau[\PB_n]}
\nc{\Mnt}{R^\tau[\M_n]}
\nc{\PBntC}{\bbC^\tau[\PB_n]}
\nc{\MntC}{\bbC^\tau[\M_n]}
\nc{\bbC}{\mathbb C}
\nc{\RSr}{R[\S_{\br}]}
\nc{\fs}{\mathfrak s}
\nc{\ft}{\mathfrak t}
\nc{\fu}{\mathfrak u}
\nc{\fv}{\mathfrak v}
\nc{\rad}{\operatorname{rad}}
\nc{\dominates}{\unrhd}

\nc{\ubluebox}[2]{\bluebox{#1}{1.7}{#2}2\udotted{#1}{#2}}
\nc{\lbluebox}[2]{\bluebox{#1}0{#2}{.3}\ldotted{#1}{#2}}
\nc{\ublueboxes}[1]{{
\foreach \x/\y in {#1}
{ \ubluebox{\x}{\y}}}
}
\nc{\lblueboxes}[1]{{
\foreach \x/\y in {#1}
{ \lbluebox{\x}{\y}}}
}

\nc{\bluebox}[4]{
\draw[color=blue!20, fill=blue!20] (#1,#2)--(#3,#2)--(#3,#4)--(#1,#4)--(#1,#2);
}
\nc{\redbox}[4]{
\draw[color=red!20, fill=red!20] (#1,#2)--(#3,#2)--(#3,#4)--(#1,#4)--(#1,#2);
}

\nc{\bluetrap}[8]{
\draw[color=blue!20, fill=blue!20] (#1,#2)--(#3,#4)--(#5,#6)--(#7,#8)--(#1,#2);
}
\nc{\redtrap}[8]{
\draw[color=red!20, fill=red!20] (#1,#2)--(#3,#4)--(#5,#6)--(#7,#8)--(#1,#2);
}

\usetikzlibrary{decorations.markings}
\usetikzlibrary{arrows,matrix}
\usepgflibrary{arrows}
\tikzset{->-/.style={decoration={
  markings,
  mark=at position #1 with {\arrow{>}}},postaction={decorate}}}
\tikzset{-<-/.style={decoration={
  markings,
  mark=at position #1 with {\arrow{<}}},postaction={decorate}}}
\nc{\Unode}[1]{\draw(#1,-2)node{$U$};}
\nc{\Dnode}[1]{\draw(#1,-2)node{$D$};}
\nc{\Fnode}[1]{\draw(#1,-2)node{$F$};}
%\nc{\Unode}[1]{\draw(#1,-2)node{$+1\phantom{+}$};}
%\nc{\Dnode}[1]{\draw(#1,-2)node{$-1\phantom{+}$};}
%\nc{\Fnode}[1]{\draw(#1,-2)node{$\phantom{+}0\phantom{+}$};}
\nc{\Cnode}[1]{\draw(#1-.1,-2)node{$\phantom{+0},$};}
\nc{\Unodes}[1]{\foreach \x in {#1}{ \Unode{\x} }}
\nc{\Dnodes}[1]{\foreach \x in {#1}{ \Dnode{\x} }}
\nc{\Fnodes}[1]{\foreach \x in {#1}{ \Fnode{\x} }}
\nc{\Cnodes}[1]{\foreach \x in {#1}{ \Cnode{\x} }}
\nc{\Uedge}[2]{\draw[->-=0.6,line width=.3mm](#1,#2-9)--(#1+1,#2+1-9); \vertsm{#1}{#2-9} \vertsm{#1+1}{#2+1-9}}
\nc{\Dedge}[2]{\draw[->-=0.6,line width=.3mm](#1,#2-9)--(#1+1,#2-1-9); \vertsm{#1}{#2-9} \vertsm{#1+1}{#2-1-9}}
\nc{\Fedge}[2]{\draw[->-=0.6,line width=.3mm](#1,#2-9)--(#1+1,#2-9); \vertsm{#1}{#2-9} \vertsm{#1+1}{#2-9}}
\nc{\Uedges}[1]{\foreach \x/\y in {#1}{\Uedge{\x}{\y}}}
\nc{\Dedges}[1]{\foreach \x/\y in {#1}{\Dedge{\x}{\y}}}
\nc{\Fedges}[1]{\foreach \x/\y in {#1}{\Fedge{\x}{\y}}}
\nc{\xvertlabel}[1]{\draw(#1,-10+.6)node{{\tiny $#1$}};}
\nc{\yvertlabel}[1]{\draw(0-.4,-9+#1)node{{\tiny $#1$}};}
\nc{\xvertlabels}[1]{\foreach \x in {#1}{ \xvertlabel{\x} }}
\nc{\yvertlabels}[1]{\foreach \x in {#1}{ \yvertlabel{\x} }}

\nc{\bbE}{\mathbb E}
\nc{\floorn}{\lfloor\tfrac n2\rfloor}
\rnc{\sp}{\supseteq}
\rnc{\arraystretch}{1.2}

\nc{\bn}{\mathbf{n}} \nc{\bt}{\mathbf{t}} \nc{\ba}{\mathbf{a}} \nc{\bl}{\mathbf{l}} \nc{\bm}{\mathbf{m}} \nc{\bk}{\mathbf{k}} \nc{\br}{\mathbf{r}} \nc{\bs}{{\mathbf s}} \nc{\bnf}{\bnf}
%\nc{\bn}{{[n]}} \nc{\bt}{{[t]}} \nc{\ba}{{[a]}} \nc{\bl}{{[l]}} \nc{\bm}{{[m]}} \nc{\bk}{{[k]}} \nc{\br}{{[r]}} \nc{\bs}{{[s]}} \nc{\bnf}{{[n-1]}}
\nc{\bp}{\mathbf{p}} 
\nc{\bq}{\mathbf{q}} 

\nc{\M}{\mathcal M}
\nc{\G}{\mathcal G}
\nc{\F}{\mathfrak F}
\nc{\MnJ}{\mathcal M_n^J}
\nc{\EnJ}{\mathcal E_n^J}
\nc{\Mat}{\operatorname{Mat}}
\nc{\RegMnJ}{\Reg(\MnJ)}
\nc{\row}{\mathfrak r}
\nc{\col}{\mathfrak c}
\nc{\Row}{\operatorname{Row}}
\nc{\Col}{\operatorname{Col}}
\nc{\Span}{\operatorname{span}}
\nc{\mat}[4]{\left[\begin{matrix}#1&#2\\#3&#4\end{matrix}\right]}
\nc{\tmat}[4]{\left[\begin{smallmatrix}#1&#2\\#3&#4\end{smallmatrix}\right]}
\nc{\ttmat}[4]{{\tiny \left[\begin{smallmatrix}#1&#2\\#3&#4\end{smallmatrix}\right]}}
\nc{\tmatt}[9]{\left[\begin{smallmatrix}#1&#2&#3\\#4&#5&#6\\#7&#8&#9\end{smallmatrix}\right]}
\nc{\ttmatt}[9]{{\tiny \left[\begin{smallmatrix}#1&#2&#3\\#4&#5&#6\\#7&#8&#9\end{smallmatrix}\right]}}
\nc{\MnGn}{\M_n\sm\G_n}
\nc{\MrGr}{\M_r\sm\G_r}
\nc{\qbin}[2]{\left[\begin{matrix}#1\\#2\end{matrix}\right]_q}
\nc{\tqbin}[2]{\left[\begin{smallmatrix}#1\\#2\end{smallmatrix}\right]_q}
\nc{\qbinx}[3]{\left[\begin{matrix}#1\\#2\end{matrix}\right]_{#3}}
\nc{\tqbinx}[3]{\left[\begin{smallmatrix}#1\\#2\end{smallmatrix}\right]_{#3}}
\nc{\MNJ}{\M_nJ}
\nc{\JMN}{J\M_n}
\nc{\RegMNJ}{\Reg(\MNJ)}
\nc{\RegJMN}{\Reg(\JMN)}
\nc{\RegMMNJ}{\Reg(\MMNJ)}
\nc{\RegJMMN}{\Reg(\JMMN)}
\nc{\Wb}{\overline{W}}
\nc{\Xb}{\overline{X}}
\nc{\Yb}{\overline{Y}}
\nc{\Zb}{\overline{Z}}
\nc{\Sib}{\overline{\Si}}
\nc{\Om}{\Omega}
\nc{\Omb}{\overline{\Om}}
\nc{\Gab}{\overline{\Ga}}
\nc{\qfact}[1]{[#1]_q!}
\nc{\smat}[2]{\left[\begin{matrix}#1&#2\end{matrix}\right]}
\nc{\tsmat}[2]{\left[\begin{smallmatrix}#1&#2\end{smallmatrix}\right]}
\nc{\hmat}[2]{\left[\begin{matrix}#1\\#2\end{matrix}\right]}
\nc{\thmat}[2]{\left[\begin{smallmatrix}#1\\#2\end{smallmatrix}\right]}
\nc{\LVW}{\mathcal L(V,W)}
\nc{\KVW}{\mathcal K(V,W)}
\nc{\LV}{\mathcal L(V)}
\nc{\RegLVW}{\Reg(\LVW)}
\nc{\sM}{\mathscr M}
\nc{\sN}{\mathscr N}
\rnc{\iff}{\ \Leftrightarrow\ }
\nc{\Hom}{\operatorname{Hom}}
\nc{\End}{\operatorname{End}}
\nc{\Aut}{\operatorname{Aut}}
\nc{\Lin}{\mathcal L}
\nc{\Hommn}{\Hom(V_m,V_n)}
\nc{\Homnm}{\Hom(V_n,V_m)}
\nc{\Homnl}{\Hom(V_n,V_l)}
\nc{\Homkm}{\Hom(V_k,V_m)}
\nc{\Endm}{\End(V_m)}
\nc{\Endn}{\End(V_n)}
\nc{\Endr}{\End(V_r)}
\nc{\Autm}{\Aut(V_m)}
\nc{\Autn}{\Aut(V_n)}
\nc{\MmnJ}{\M_{mn}^J}
\nc{\MmnA}{\M_{mn}^A}
\nc{\MmnB}{\M_{mn}^B}
\nc{\Mmn}{\M_{mn}}
\nc{\Mkl}{\M_{kl}}
\nc{\Mnm}{\M_{nm}}
\nc{\EmnJ}{\mathcal E_{mn}^J}
\nc{\MmGm}{\M_m\sm\G_m}
\nc{\RegMmnJ}{\Reg(\MmnJ)}
\rnc{\implies}{\ \Rightarrow\ }
\nc{\DMmn}[1]{D_{#1}(\Mmn)}
\nc{\DMmnJ}[1]{D_{#1}(\MmnJ)}
\nc{\MMNJ}{\Mmn J}
\nc{\JMMN}{J\Mmn}
\nc{\JMMNJ}{J\Mmn J}
\nc{\Inr}{\mathcal I(V_n,W_r)}
\nc{\Lnr}{\mathcal L(V_n,W_r)}
\nc{\Knr}{\mathcal K(V_n,W_r)}
\nc{\Imr}{\mathcal I(V_m,W_r)}
\nc{\Kmr}{\mathcal K(V_m,W_r)}
\nc{\Lmr}{\mathcal L(V_m,W_r)}
\nc{\Kmmr}{\mathcal K(V_m,W_{m-r})}
\nc{\tr}{{\operatorname{T}}}
\nc{\MMN}{\MmnA(\F_1)}
\nc{\MKL}{\Mkl^B(\F_2)}
\nc{\RegMMN}{\Reg(\MmnA(\F_1))}
\nc{\RegMKL}{\Reg(\Mkl^B(\F_2))}
\nc{\gRhA}{\widehat{\mathscr R}^A}
\nc{\gRhB}{\widehat{\mathscr R}^B}
\nc{\gLhA}{\widehat{\mathscr L}^A}
\nc{\gLhB}{\widehat{\mathscr L}^B}
\nc{\timplies}{\Rightarrow}
\nc{\tiff}{\Leftrightarrow}
\nc{\Sija}{S_{ij}^a}
\nc{\dmat}[8]{\draw(#1*1.5,#2)node{$\left[\begin{smallmatrix}#3&#4&#5\\#6&#7&#8\end{smallmatrix}\right]$};}
\nc{\bdmat}[8]{\draw(#1*1.5,#2)node{${\mathbf{\left[\begin{smallmatrix}#3&#4&#5\\#6&#7&#8\end{smallmatrix}\right]}}$};}
\nc{\rdmat}[8]{\draw(#1*1.5,#2)node{\rotatebox{90}{$\left[\begin{smallmatrix}#3&#4&#5\\#6&#7&#8\end{smallmatrix}\right]$}};}
\nc{\rldmat}[8]{\draw(#1*1.5-0.375,#2)node{\rotatebox{90}{$\left[\begin{smallmatrix}#3&#4&#5\\#6&#7&#8\end{smallmatrix}\right]$}};}
\nc{\rrdmat}[8]{\draw(#1*1.5+.375,#2)node{\rotatebox{90}{$\left[\begin{smallmatrix}#3&#4&#5\\#6&#7&#8\end{smallmatrix}\right]$}};}
\nc{\rfldmat}[8]{\draw(#1*1.5-0.375+.15,#2)node{\rotatebox{90}{$\left[\begin{smallmatrix}#3&#4&#5\\#6&#7&#8\end{smallmatrix}\right]$}};}
\nc{\rfrdmat}[8]{\draw(#1*1.5+.375-.15,#2)node{\rotatebox{90}{$\left[\begin{smallmatrix}#3&#4&#5\\#6&#7&#8\end{smallmatrix}\right]$}};}
\nc{\xL}{[x]_{\! _\gL}}\nc{\yL}{[y]_{\! _\gL}}\nc{\xR}{[x]_{\! _\gR}}\nc{\yR}{[y]_{\! _\gR}}\nc{\xH}{[x]_{\! _\gH}}\nc{\yH}{[y]_{\! _\gH}}\nc{\XK}{[X]_{\! _\gK}}\nc{\xK}{[x]_{\! _\gK}}
%\nc{\xL}{[x]_\L}\nc{\yL}{[y]_\L}\nc{\xR}{[x]_\R}\nc{\yR}{[y]_\R}\nc{\xH}{[x]_\H}\nc{\yH}{[y]_\H}
\nc{\RegSija}{\Reg(\Sija)}
\nc{\MnmK}{\M_{nm}^K}
\nc{\cC}{\mathcal C}
\nc{\cR}{\mathcal R}
\nc{\Ckl}{\cC_k(l)}
\nc{\Rkl}{\cR_k(l)}
\nc{\Cmr}{\cC_m(r)}
\nc{\Rmr}{\cR_m(r)}
\nc{\Cnr}{\cC_n(r)}
\nc{\Rnr}{\cR_n(r)}
\nc{\Z}{\mathcal Z}

\nc{\Reg}{\operatorname{Reg}}
%\nc{\RP}{\operatorname{RP}}
\nc{\TXa}{\T_X^a}
\nc{\TXA}{\T(X,A)}
\nc{\TXal}{\T(X,\al)}
\nc{\RegTXa}{\Reg(\TXa)}
\nc{\RegTXA}{\Reg(\TXA)}
\nc{\RegTXal}{\Reg(\TXal)}
\nc{\PalX}{\P_\al(X)}
\nc{\EAX}{\E_A(X)}
\nc{\Bb}{\overline{B}}
\nc{\bb}{\overline{\be}}
\nc{\bw}{{\bf w}}
\nc{\bz}{{\bf z}}
\nc{\TASA}{\T_A\sm\S_A}
\nc{\Ub}{\overline{U}}
\nc{\Vb}{\overline{V}}
\nc{\eb}{\overline{e}}
\nc{\ob}{\overline{o}}
\nc{\cb}{\overline{c}}
\nc{\db}{\overline{d}}
\nc{\tb}{\overline{t}}
\nc{\Eb}{\overline{E}}
\nc{\Ob}{\overline{O}}
\nc{\Qb}{\overline{Q}}
\nc{\EXa}{\E_X^a}
\nc{\oijr}{1\leq i<j\leq r}
\nc{\veb}{\overline{\ve}}
\nc{\bbT}{\mathbb T}
\nc{\Surj}{\operatorname{Surj}}
\nc{\Sone}{S^{(1)}}
\nc{\fillbox}[2]{\draw[fill=gray!30](#1,#2)--(#1+1,#2)--(#1+1,#2+1)--(#1,#2+1)--(#1,#2);}
\nc{\raa}{\rangle_J}
\nc{\raJ}{\rangle_J}
%\rnc{\star}{\star}
\nc{\Ea}{E_J}
\nc{\EJ}{E_J}
\nc{\ep}{\epsilon} \nc{\ve}{\varepsilon}
%\nc{\ep}{\varepsilon} \nc{\ve}{\eta}
\nc{\IXa}{\I_X^a}
\nc{\RegIXa}{\Reg(\IXa)}
\nc{\JXa}{\J_X^a}
\nc{\RegJXa}{\Reg(\JXa)}
\nc{\IXA}{\I(X,A)}
\nc{\IAX}{\I(A,X)}
\nc{\RegIXA}{\Reg(\IXA)}
\nc{\RegIAX}{\Reg(\IAX)}
\nc{\trans}[2]{\left(\begin{smallmatrix} #1 \\ #2 \end{smallmatrix}\right)}
\nc{\bigtrans}[2]{\left(\begin{matrix} #1 \\ #2 \end{matrix}\right)}
\nc{\lmap}[1]{\mapstochar \xrightarrow {\ #1\ }}
\nc{\EaTXa}{E}

\nc{\gL}{\mathscr L}
\nc{\gR}{\mathscr R}
\nc{\gH}{\mathscr H}
\nc{\gJ}{\mathscr J}
\nc{\gD}{\mathscr D}
\nc{\gK}{\mathscr K}
\nc{\gLa}{\mathscr L^a}
\nc{\gRa}{\mathscr R^a}
\nc{\gHa}{\mathscr H^a}
\nc{\gJa}{\mathscr J^a}
\nc{\gDa}{\mathscr D^a}
\nc{\gKa}{\mathscr K^a}
\nc{\gLJ}{\mathscr L^J}
\nc{\gRJ}{\mathscr R^J}
\nc{\gHJ}{\mathscr H^J}
\nc{\gJJ}{\mathscr J^J}
\nc{\gDJ}{\mathscr D^J}
\nc{\gKJ}{\mathscr K^J}
\nc{\gLh}{\widehat{\mathscr L}^J}
\nc{\gRh}{\widehat{\mathscr R}^J}
\nc{\gHh}{\widehat{\mathscr H}^J}
\nc{\gJh}{\widehat{\mathscr J}^J}
\nc{\gDh}{\widehat{\mathscr D}^J}
\nc{\gKh}{\widehat{\mathscr K}^J}
\nc{\Lh}{\widehat{L}^J}
\nc{\Rh}{\widehat{R}}
\nc{\Hh}{\widehat{H}^J}
\nc{\Jh}{\widehat{J}^J}
\nc{\Dh}{\widehat{D}^J}
\nc{\Kh}{\widehat{K}^J}
\nc{\gLb}{\widehat{\mathscr L}}
\nc{\gRb}{\widehat{\mathscr R}}
\nc{\gHb}{\widehat{\mathscr H}}
\nc{\gJb}{\widehat{\mathscr J}}
\nc{\gDb}{\widehat{\mathscr D}}
\nc{\gKb}{\widehat{\mathscr K}}
\nc{\Lb}{\widehat{L}^J}
\nc{\Rb}{\widehat{R}^J}
\nc{\Hb}{\widehat{H}^J}
\nc{\Jb}{\widehat{J}^J}
\nc{\Db}{\overline{D}}
\nc{\Kb}{\widehat{K}}

\hyphenation{mon-oid mon-oids}

\nc{\itemit}[1]{\item[\emph{(#1)}]}
\nc{\E}{\mathcal E}
\nc{\TX}{\T(X)}
\nc{\TXP}{\T(X,\P)}
\nc{\EX}{\E(X)}
\nc{\EXP}{\E(X,\P)}
\nc{\SX}{\S(X)}
\nc{\SXP}{\S(X,\P)}
\nc{\Sing}{\operatorname{Sing}}
%\nc{\Sing}{\E}
\nc{\idrank}{\operatorname{idrank}}
\nc{\SingXP}{\Sing(X,\P)}
\nc{\De}{\Delta}
\nc{\sgp}{\operatorname{sgp}}
\nc{\mon}{\operatorname{mon}}
\nc{\Dn}{\mathcal D_n}
\nc{\Dm}{\mathcal D_m}

\nc{\lline}[1]{\draw(3*#1,0)--(3*#1+2,0);}
\nc{\uline}[1]{\draw(3*#1,5)--(3*#1+2,5);}
\nc{\thickline}[2]{\draw(3*#1,5)--(3*#2,0); \draw(3*#1+2,5)--(3*#2+2,0) ;}
\nc{\thicklabel}[3]{\draw(3*#1+1+3*#2*0.15-3*#1*0.15,4.25)node{{\tiny $#3$}};}

\nc{\slline}[3]{\draw(3*#1+#3,0+#2)--(3*#1+2+#3,0+#2);}
\nc{\suline}[3]{\draw(3*#1+#3,5+#2)--(3*#1+2+#3,5+#2);}
\nc{\sthickline}[4]{\draw(3*#1+#4,5+#3)--(3*#2+#4,0+#3); \draw(3*#1+2+#4,5+#3)--(3*#2+2+#4,0+#3) ;}
\nc{\sthicklabel}[5]{\draw(3*#1+1+3*#2*0.15-3*#1*0.15+#5,4.25+#4)node{{\tiny $#3$}};}

\nc{\stll}[5]{\sthickline{#1}{#2}{#4}{#5} \sthicklabel{#1}{#2}{#3}{#4}{#5}}
\nc{\tll}[3]{\stll{#1}{#2}{#3}00}

\nc{\mfourpic}[9]{
\slline1{#9}0
\slline3{#9}0
\slline4{#9}0
\slline5{#9}0
\suline1{#9}0
\suline3{#9}0
\suline4{#9}0
\suline5{#9}0
\stll1{#1}{#5}{#9}{0}
\stll3{#2}{#6}{#9}{0}
\stll4{#3}{#7}{#9}{0}
\stll5{#4}{#8}{#9}{0}
\draw[dotted](6,0+#9)--(8,0+#9);
\draw[dotted](6,5+#9)--(8,5+#9);
}
\nc{\vdotted}[1]{
\draw[dotted](3*#1,10)--(3*#1,15);
\draw[dotted](3*#1+2,10)--(3*#1+2,15);
}

\nc{\Clab}[2]{
\sthicklabel{#1}{#1}{{}_{\phantom{#1}}C_{#1}}{1.25+5*#2}0
}
\nc{\sClab}[3]{
\sthicklabel{#1}{#1}{{}_{\phantom{#1}}C_{#1}}{1.25+5*#2}{#3}
}
\nc{\Clabl}[3]{
\sthicklabel{#1}{#1}{{}_{\phantom{#3}}C_{#3}}{1.25+5*#2}0
}
\nc{\sClabl}[4]{
\sthicklabel{#1}{#1}{{}_{\phantom{#4}}C_{#4}}{1.25+5*#2}{#3}
}
\nc{\Clabll}[3]{
\sthicklabel{#1}{#1}{C_{#3}}{1.25+5*#2}0
}
\nc{\sClabll}[4]{
\sthicklabel{#1}{#1}{C_{#3}}{1.25+5*#2}{#3}
}

\nc{\mtwopic}[6]{
\slline1{#6*5}{#5}
\slline2{#6*5}{#5}
\suline1{#6*5}{#5}
\suline2{#6*5}{#5}
\stll1{#1}{#3}{#6*5}{#5}
\stll2{#2}{#4}{#6*5}{#5}
}
\nc{\mtwopicl}[6]{
\slline1{#6*5}{#5}
\slline2{#6*5}{#5}
\suline1{#6*5}{#5}
\suline2{#6*5}{#5}
\stll1{#1}{#3}{#6*5}{#5}
\stll2{#2}{#4}{#6*5}{#5}
\sClabl1{#6}{#5}{i}
\sClabl2{#6}{#5}{j}
}

%\nc{\keru}{\ker_{\operatorname{u}}} \nc{\kerl}{\ker_{\operatorname{l}}}
%\nc{\dom}{\operatorname{dom}_u} \nc{\codom}{\operatorname{dom}_l}
\nc{\keru}{\operatorname{ker}^\wedge} \nc{\kerl}{\operatorname{ker}_\vee}%\nc{\codom}{\operatorname{im}}

\nc{\coker}{\operatorname{coker}}
%\nc{\KER}{\operatorname{KER}}
\nc{\KER}{\ker}
\nc{\N}{\mathbb N}
\nc{\LaBn}{L_\al(\B_n)}
\nc{\RaBn}{R_\al(\B_n)}
\nc{\LaPBn}{L_\al(\PB_n)}
\nc{\RaPBn}{R_\al(\PB_n)}
\nc{\rhorBn}{\rho_r(\B_n)}
\nc{\DrBn}{D_r(\B_n)}
\nc{\DrPn}{D_r(\P_n)}
\nc{\DrPBn}{D_r(\PB_n)}
\nc{\DrKn}{D_r(\K_n)}
\nc{\alb}{\al_{\vee}}
\nc{\beb}{\be^{\wedge}}
\nc{\Bal}{\operatorname{Bal}}
\nc{\Red}{\operatorname{Red}}
\nc{\Pnxi}{\P_n^\xi}
\nc{\Bnxi}{\B_n^\xi}
\nc{\PBnxi}{\PB_n^\xi}
\nc{\Knxi}{\K_n^\xi}
\nc{\C}{\mathbb C}
\nc{\exi}{e^\xi}
\nc{\Exi}{E^\xi}
\nc{\eximu}{e^\xi_\mu}
\nc{\Eximu}{E^\xi_\mu}
\nc{\REF}{ {\red [Ref?]} }
\nc{\GL}{\operatorname{GL}}
%\rnc{\O}{\operatorname{O}}
\rnc{\O}{\mathcal O}

\nc{\vtx}[2]{\fill (#1,#2)circle(.2);}
\nc{\lvtx}[2]{\fill (#1,0)circle(.2);}
\nc{\uvtx}[2]{\fill (#1,1.5)circle(.2);}

\nc{\Eq}{\mathfrak{Eq}}
%\nc{\Gau}{\Ga_{\operatorname{u}}} \nc{\Gal}{\Ga_{\operatorname{l}}}
\nc{\Gau}{\Ga^\wedge} \nc{\Gal}{\Ga_\vee}
%\nc{\Lamu}{\Lam_{\operatorname{u}}} \nc{\Laml}{\Lam_{\operatorname{l}}}
\nc{\Lamu}{\Lam^\wedge} \nc{\Laml}{\Lam_\vee}
\nc{\bX}{{\bf X}}
\nc{\bY}{{\bf Y}}
\nc{\ds}{\displaystyle}

\nc{\uuvert}[1]{\fill (#1,3)circle(.2);}
\nc{\uuuvert}[1]{\fill (#1,4.5)circle(.2);}
\nc{\overt}[1]{\fill (#1,0)circle(.1);}
\nc{\overtl}[3]{\node[vertex] (#3) at (#1,0) {  {\tiny $#2$} };}
\nc{\cv}[2]{\draw(#1,1.5) to [out=270,in=90] (#2,0);}
\nc{\cvs}[2]{\draw(#1,1.5) to [out=270+30,in=90+30] (#2,0);}
\nc{\ucv}[2]{\draw(#1,3) to [out=270,in=90] (#2,1.5);}
\nc{\uucv}[2]{\draw(#1,4.5) to [out=270,in=90] (#2,3);}
\nc{\textpartn}[1]{{\lower1.0 ex\hbox{\begin{tikzpicture}[xscale=.3,yscale=0.3] #1 \end{tikzpicture}}}}
\nc{\textpartnx}[2]{{\lower1.0 ex\hbox{\begin{tikzpicture}[xscale=.3,yscale=0.3] 
\foreach \x in {1,...,#1}
{ \uvert{\x} \lvert{\x} }
#2 \end{tikzpicture}}}}
\nc{\disppartnx}[2]{{\lower1.0 ex\hbox{\begin{tikzpicture}[scale=0.3] 
\foreach \x in {1,...,#1}
{ \uvert{\x} \lvert{\x} }
#2 \end{tikzpicture}}}}
\nc{\disppartnxd}[2]{{\lower2.1 ex\hbox{\begin{tikzpicture}[scale=0.3] 
\foreach \x in {1,...,#1}
{ \uuvert{\x} \uvert{\x} \lvert{\x} }
#2 \end{tikzpicture}}}}
\nc{\disppartnxdn}[2]{{\lower2.1 ex\hbox{\begin{tikzpicture}[scale=0.3] 
\foreach \x in {1,...,#1}
{ \uuvert{\x} \lvert{\x} }
#2 \end{tikzpicture}}}}
%\nc{\disppartnxdd}[2]{{\lower2.1 ex\hbox{\begin{tikzpicture}[scale=0.3] 
%\foreach \x in {1,...,#1}
%{ \uuuvert{\x} \uuvert{\x} \uvert{\x} \lvert{\x} }
%#2 \end{tikzpicture}}}}
\nc{\disppartnxdd}[2]{{\lower3.6 ex\hbox{\begin{tikzpicture}[scale=0.3] 
\foreach \x in {1,...,#1}
{ \uuuvert{\x} \uuvert{\x} \uvert{\x} \lvert{\x} }
#2 \end{tikzpicture}}}}

\nc{\dispgax}[2]{{\lower0.0 ex\hbox{\begin{tikzpicture}[scale=0.3] 
#2
\foreach \x in {1,...,#1}
{\lvert{\x} }
 \end{tikzpicture}}}}
\nc{\textgax}[2]{{\lower0.4 ex\hbox{\begin{tikzpicture}[scale=0.3] 
#2
\foreach \x in {1,...,#1}
{\lvert{\x} }
 \end{tikzpicture}}}}
\nc{\textlinegraph}[2]{{\raise#1 ex\hbox{\begin{tikzpicture}[scale=0.8] 
#2
 \end{tikzpicture}}}}
\nc{\textlinegraphl}[2]{{\raise#1 ex\hbox{\begin{tikzpicture}[scale=0.8] 
\tikzstyle{vertex}=[circle,draw=black, fill=white, inner sep = 0.07cm]
#2
 \end{tikzpicture}}}}
\nc{\displinegraph}[1]{{\lower0.0 ex\hbox{\begin{tikzpicture}[scale=0.6] 
#1
 \end{tikzpicture}}}}
 
\nc{\disppartnthreeone}[1]{{\lower1.0 ex\hbox{\begin{tikzpicture}[scale=0.3] 
\foreach \x in {1,2,3,5,6}
{ \uvert{\x} }
\foreach \x in {1,2,4,5,6}
{ \lvert{\x} }
\draw[dotted] (3.5,1.5)--(4.5,1.5);
\draw[dotted] (2.5,0)--(3.5,0);
#1 \end{tikzpicture}}}}

\nc{\partn}[4]{\left( \begin{array}{c|c} %fine
#1 \ & \ #3 \ \ \\ \cline{2-2}
#2 \ & \ #4 \ \
\end{array} \!\!\! \right)}
\nc{\partnlong}[6]{\partn{#1}{#2}{#3,\ #4}{#5,\ #6}} %fine
\nc{\partnsh}[2]{\left( \begin{array}{c} %fine
#1 \\
#2 
\end{array} \right)}
\nc{\partncodefz}[3]{\partn{#1}{#2}{#3}{\emptyset}}
\nc{\partndefz}[3]{{\partn{#1}{#2}{\emptyset}{#3}}}
\nc{\partnlast}[2]{\left( \begin{array}{c|c}
#1 \ &  \ #2 \\
#1 \ &  \ #2
\end{array} \right)}

\nc{\partnlist}[8]{
\left( \begin{array}{c|c|c|c|c|c} %fine
\!\! #1 & \cdots & #2 & #3 & \cdots & #4\ \ \\ \cline{4-6}
\!\! #5 & \cdots & #6 & #7 & \cdots & #8 \ \
\end{array} \!\!\! \right)
}

\nc{\partnlistnodef}[8]{
\left( \begin{array}{c|c|c|c|c|c} %fine
\!\! #1 & \cdots & #2 & #3 & \cdots & #4\ \ \\ %\cline{4-6}
\!\! #5 & \cdots & #6 & #7 & \cdots & #8 \ \
\end{array} \!\!\! \right)
}

\nc{\partnlistnodefshort}[4]{
\left( \begin{array}{c|c|c} %fine
\!\! #1 & \cdots & #2 \\ 
\!\! #3 & \cdots & #4 
\end{array} \!\!\! \right)
}

\nc{\partialpermshort}[4]{
\left[ \begin{array}{c|c|c} %fine
\! #1 & \cdots & #2 \\ 
\! #3 & \cdots & #4 
\end{array} \!\! \right]
}

\nc{\partnlistnodefshortish}[6]{
\left( \begin{array}{c|c|c|c} %fine
\!\! #1 & \cdots & #2 & #3 \\ 
\!\! #4 & \cdots & #5 & #6 
\end{array} \!\!\! \right)
}

\nc{\partialpermshortish}[6]{
\left[ \begin{array}{c|c|c|c} %fine
\!\! #1 & \cdots & #2 & #3 \\ 
\!\! #4 & \cdots & #5 & #6 
\end{array} \!\!\! \right]
}

\nc{\uuarcx}[3]{\draw(#1,3)arc(180:270:#3) (#1+#3,3-#3)--(#2-#3,3-#3) (#2-#3,3-#3) arc(270:360:#3);}
\nc{\uuarc}[2]{\uuarcx{#1}{#2}{.4}}
\nc{\uuuarcx}[3]{\draw(#1,4.5)arc(180:270:#3) (#1+#3,4.5-#3)--(#2-#3,4.5-#3) (#2-#3,4.5-#3) arc(270:360:#3);}
\nc{\uuuarc}[2]{\uuuarcx{#1}{#2}{.4}}
%\nc{\darcx}[3]{\draw(#1,0)arc(180:90:#3) (#1+#3,#3)--(#2-#3,#3) (#2-#3,#3) arc(90:0:#3);}
%\nc{\darc}[2]{\darcx{#1}{#2}{.4}}
\nc{\udarcx}[3]{\draw(#1,1.5)arc(180:90:#3) (#1+#3,1.5+#3)--(#2-#3,1.5+#3) (#2-#3,1.5+#3) arc(90:0:#3);}
\nc{\udarc}[2]{\udarcx{#1}{#2}{.4}}
\nc{\uudarcx}[3]{\draw(#1,3)arc(180:90:#3) (#1+#3,3+#3)--(#2-#3,3+#3) (#2-#3,3+#3) arc(90:0:#3);}
\nc{\uudarc}[2]{\uudarcx{#1}{#2}{.4}}
%\nc{\uarcx}[3]{\draw(#1,1.5)arc(180:270:#3) (#1+#3,1.5-#3)--(#2-#3,1.5-#3) (#2-#3,1.5-#3) arc(270:360:#3);}
%\nc{\uarc}[2]{\uarcx{#1}{#2}{.4}}
\nc{\darcxhalf}[3]{\draw(#1,0)arc(180:90:#3) (#1+#3,#3)--(#2,#3) ;}
\nc{\darchalf}[2]{\darcxhalf{#1}{#2}{.4}}
\nc{\uarcxhalf}[3]{\draw(#1,1.5)arc(180:270:#3) (#1+#3,1.5-#3)--(#2,1.5-#3) ;}
\nc{\uarchalf}[2]{\uarcxhalf{#1}{#2}{.4}}
\nc{\uarcxhalfr}[3]{\draw (#1+#3,1.5-#3)--(#2-#3,1.5-#3) (#2-#3,1.5-#3) arc(270:360:#3);}
\nc{\uarchalfr}[2]{\uarcxhalfr{#1}{#2}{.4}}

\nc{\bdarcx}[3]{\draw[blue](#1,0)arc(180:90:#3) (#1+#3,#3)--(#2-#3,#3) (#2-#3,#3) arc(90:0:#3);}
\nc{\bdarc}[2]{\darcx{#1}{#2}{.4}}
\nc{\rduarcx}[3]{\draw[red](#1,0)arc(180:270:#3) (#1+#3,0-#3)--(#2-#3,0-#3) (#2-#3,0-#3) arc(270:360:#3);}
\nc{\rduarc}[2]{\uarcx{#1}{#2}{.4}}
\nc{\duarcx}[3]{\draw(#1,0)arc(180:270:#3) (#1+#3,0-#3)--(#2-#3,0-#3) (#2-#3,0-#3) arc(270:360:#3);}
\nc{\duarc}[2]{\uarcx{#1}{#2}{.4}}

\nc{\uuv}[1]{\fill (#1,4)circle(.1);}
\nc{\uv}[1]{\fill (#1,2)circle(.1);}
\nc{\lv}[1]{\fill (#1,0)circle(.1);}
\nc{\uvw}[1]{\draw[fill=white] (#1,2)circle(.18);}
\nc{\lvw}[1]{\draw (#1,0)circle(.18);}
\nc{\uvred}[1]{\fill[red] (#1,2)circle(.1);}
\nc{\lvred}[1]{\fill[red] (#1,0)circle(.1);}
\nc{\lvwhite}[1]{\fill[white] (#1,0)circle(.1);}
\nc{\buv}[1]{\fill (#1,2)circle(.18);}
\nc{\blv}[1]{\fill (#1,0)circle(.18);}

\nc{\uvs}[1]{{
\foreach \x in {#1}
{ \uv{\x}}
}}
\nc{\uuvs}[1]{{
\foreach \x in {#1}
{ \uuv{\x}}
}}
\nc{\lvs}[1]{{
\foreach \x in {#1}
{ \lv{\x}}
}}

\nc{\buvs}[1]{{
\foreach \x in {#1}
{ \buv{\x}}
}}
\nc{\blvs}[1]{{
\foreach \x in {#1}
{ \blv{\x}}
}}

\nc{\uvreds}[1]{{
\foreach \x in {#1}
{ \uvred{\x}}
}}
\nc{\lvreds}[1]{{
\foreach \x in {#1}
{ \lvred{\x}}
}}

\nc{\uudotted}[2]{\draw [dotted] (#1,4)--(#2,4);}
\nc{\uudotteds}[1]{{
\foreach \x/\y in {#1}
{ \uudotted{\x}{\y}}
}}
\nc{\uudottedsm}[2]{\draw [dotted] (#1+.4,4)--(#2-.4,4);}
\nc{\uudottedsms}[1]{{
\foreach \x/\y in {#1}
{ \uudottedsm{\x}{\y}}
}}
\nc{\udottedsm}[2]{\draw [dotted] (#1+.4,2)--(#2-.4,2);}
\nc{\udottedsms}[1]{{
\foreach \x/\y in {#1}
{ \udottedsm{\x}{\y}}
}}
\nc{\udotted}[2]{\draw [dotted] (#1,2)--(#2,2);}
\nc{\udotteds}[1]{{
\foreach \x/\y in {#1}
{ \udotted{\x}{\y}}
}}
\nc{\ldotted}[2]{\draw [dotted] (#1,0)--(#2,0);}
\nc{\ldotteds}[1]{{
\foreach \x/\y in {#1}
{ \ldotted{\x}{\y}}
}}
\nc{\ldottedsm}[2]{\draw [dotted] (#1+.4,0)--(#2-.4,0);}
\nc{\ldottedsms}[1]{{
\foreach \x/\y in {#1}
{ \ldottedsm{\x}{\y}}
}}

\nc{\stlinest}[2]{\draw(#1,4)--(#2,0);}

\nc{\stlined}[2]{\draw[dotted](#1,2)--(#2,0);}

%\nc{\stline}[2]{\draw(#1,2)--(#2,0);}
\nc{\tlab}[2]{\draw(#1,2)node[above]{\tiny $#2$};}
\nc{\tudots}[1]{\draw(#1,2)node{$\cdots$};}
\nc{\tldots}[1]{\draw(#1,0)node{$\cdots$};}

%\nc{\uvw}[1]{\fill[white] (#1,2)circle(.1);}
\nc{\huv}[1]{\fill (#1,1)circle(.1);}
\nc{\llv}[1]{\fill (#1,-2)circle(.1);}
\nc{\arcup}[2]{
\draw(#1,2)arc(180:270:.4) (#1+.4,1.6)--(#2-.4,1.6) (#2-.4,1.6) arc(270:360:.4);
}
\nc{\harcup}[2]{
\draw(#1,1)arc(180:270:.4) (#1+.4,.6)--(#2-.4,.6) (#2-.4,.6) arc(270:360:.4);
}
\nc{\arcdn}[2]{
\draw(#1,0)arc(180:90:.4) (#1+.4,.4)--(#2-.4,.4) (#2-.4,.4) arc(90:0:.4);
}
\nc{\cve}[2]{
\draw(#1,2) to [out=270,in=90] (#2,0);
}
\nc{\hcve}[2]{
\draw(#1,1) to [out=270,in=90] (#2,0);
}
\nc{\catarc}[3]{
\draw(#1,2)arc(180:270:#3) (#1+#3,2-#3)--(#2-#3,2-#3) (#2-#3,2-#3) arc(270:360:#3);
}

\nc{\arcr}[2]{
\draw[red](#1,0)arc(180:90:.4) (#1+.4,.4)--(#2-.4,.4) (#2-.4,.4) arc(90:0:.4);
}
\nc{\arcb}[2]{
\draw[blue](#1,2-2)arc(180:270:.4) (#1+.4,1.6-2)--(#2-.4,1.6-2) (#2-.4,1.6-2) arc(270:360:.4);
}
\nc{\loopr}[1]{
\draw[blue](#1,-2) edge [out=130,in=50,loop] ();
}
\nc{\loopb}[1]{
\draw[red](#1,-2) edge [out=180+130,in=180+50,loop] ();
}
%\nc{\arcr}[2]{
%\draw[red](#1,0-2)arc(180:90:.4) (#1+.4,.4-2)--(#2-.4,.4-2) (#2-.4,.4-2) arc(90:0:.4);
%}
%\nc{\arcb}[2]{
%\draw[blue](#1,2-2-2)arc(180:270:.4) (#1+.4,1.6-2-2)--(#2-.4,1.6-2-2) (#2-.4,1.6-2-2) arc(270:360:.4);
%}
%\nc{\loopr}[1]{
%\draw[red](#1,0-2) edge [out=130,in=50,loop] ();
%}
%\nc{\loopb}[1]{
%\draw[blue](#1,0-2) edge [out=180+130,in=180+50,loop] ();
%}
\nc{\redto}[2]{\draw[red](#1,0)--(#2,0);}
\nc{\bluto}[2]{\draw[blue](#1,0)--(#2,0);}
\nc{\dotto}[2]{\draw[dotted](#1,0)--(#2,0);}
\nc{\lloopr}[2]{\draw[red](#1,0)arc(0:360:#2);}
\nc{\lloopb}[2]{\draw[blue](#1,0)arc(0:360:#2);}
\nc{\rloopr}[2]{\draw[red](#1,0)arc(-180:180:#2);}
\nc{\rloopb}[2]{\draw[blue](#1,0)arc(-180:180:#2);}
\nc{\uloopr}[2]{\draw[red](#1,0)arc(-270:270:#2);}
\nc{\uloopb}[2]{\draw[blue](#1,0)arc(-270:270:#2);}
\nc{\dloopr}[2]{\draw[red](#1,0)arc(-90:270:#2);}
\nc{\dloopb}[2]{\draw[blue](#1,0)arc(-90:270:#2);}
\nc{\llloopr}[2]{\draw[red](#1,0-2)arc(0:360:#2);}
\nc{\llloopb}[2]{\draw[blue](#1,0-2)arc(0:360:#2);}
\nc{\lrloopr}[2]{\draw[red](#1,0-2)arc(-180:180:#2);}
\nc{\lrloopb}[2]{\draw[blue](#1,0-2)arc(-180:180:#2);}
\nc{\ldloopr}[2]{\draw[red](#1,0-2)arc(-270:270:#2);}
\nc{\ldloopb}[2]{\draw[blue](#1,0-2)arc(-270:270:#2);}
\nc{\luloopr}[2]{\draw[red](#1,0-2)arc(-90:270:#2);}
\nc{\luloopb}[2]{\draw[blue](#1,0-2)arc(-90:270:#2);}

\nc{\larcb}[2]{
\draw[blue](#1,0-2)arc(180:90:.4) (#1+.4,.4-2)--(#2-.4,.4-2) (#2-.4,.4-2) arc(90:0:.4);
}
\nc{\larcr}[2]{
\draw[red](#1,2-2-2)arc(180:270:.4) (#1+.4,1.6-2-2)--(#2-.4,1.6-2-2) (#2-.4,1.6-2-2) arc(270:360:.4);
}

\rnc{\H}{\mathscr H}
\rnc{\L}{\mathscr L}
\nc{\R}{\mathcal R}
\nc{\D}{\mathscr D}
%\nc{\J}{\mathscr J}
\nc{\J}{\mathcal J}

\nc{\ssim}{\mathrel{\raise0.25 ex\hbox{\oalign{$\approx$\crcr\noalign{\kern-0.84 ex}$\approx$}}}}
\nc{\POI}{\mathcal{O}}
%\nc{\POI}{\mathcal{POI}}
\nc{\wb}{\overline{w}}
\nc{\zb}{\overline{z}}
\nc{\ub}{\overline{u}}
\nc{\vb}{\overline{v}}
\nc{\fb}{\overline{f}}
\nc{\gb}{\overline{g}}
\nc{\hb}{\overline{h}}
\nc{\pb}{\overline{p}}
\nc{\xb}{\overline{x}}
\nc{\qb}{\overline{q}}
\rnc{\sb}{\overline{s}}
\nc{\Sb}{\overline{\si}}
\nc{\XR}{\pres{X}{R\,}}
\nc{\YQ}{\pres{Y}{Q}}
\nc{\ZP}{\pres{Z}{P\,}}
\nc{\XRone}{\pres{X_1}{R_1}}
\nc{\XRtwo}{\pres{X_2}{R_2}}
\nc{\XRthree}{\pres{X_1\cup X_2}{R_1\cup R_2\cup R_3}}
\nc{\er}{\eqref}
\nc{\larr}{\mathrel{\hspace{-0.35 ex}>\hspace{-1.1ex}-}\hspace{-0.35 ex}}
\nc{\rarr}{\mathrel{\hspace{-0.35 ex}-\hspace{-0.5ex}-\hspace{-2.3ex}>\hspace{-0.35 ex}}}
\nc{\lrarr}{\mathrel{\hspace{-0.35 ex}>\hspace{-1.1ex}-\hspace{-0.5ex}-\hspace{-2.3ex}>\hspace{-0.35 ex}}}
\nc{\nn}{\tag*{}}
\nc{\epfal}{\tag*{$\Box$}}
\nc{\tagd}[1]{\tag*{(#1)$'$}}
\nc{\ldb}{[\![}
\nc{\rdb}{]\!]}
\nc{\sm}{\setminus}
\nc{\I}{\mathcal I}
\nc{\InSn}{\I_n\setminus\S_n}
%\nc{\dom}{\operatorname{dom}_{\operatorname{u}}} \nc{\codom}{\operatorname{dom}_{\operatorname{l}}}
%\nc{\dom}{\operatorname{dom}_u} \nc{\codom}{\operatorname{dom}_l}
\nc{\dom}{\operatorname{dom}} \nc{\codom}{\operatorname{codom}}%\nc{\codom}{\operatorname{im}}
\nc{\ojin}{1\leq j<i\leq n}
%\nc{\R}{\mathcal R}
%\rnc{\L}{\mathcal L}
\nc{\eh}{\widehat{e}}
\nc{\wh}{\widehat{w}}
\nc{\uh}{\widehat{u}}
\nc{\vh}{\widehat{v}}
\nc{\sh}{\widehat{s}}
\nc{\fh}{\widehat{f}}
\nc{\textres}[1]{\text{\emph{#1}}}
\nc{\aand}{\emph{\ and \ }}
\nc{\iif}{\emph{\ if \ }}
\nc{\textlarr}{\ \larr\ }
\nc{\textrarr}{\ \rarr\ }
\nc{\textlrarr}{\ \lrarr\ }

\nc{\comma}{,\ }

\nc{\COMMA}{,\quad}
\nc{\TnSn}{\T_n\setminus\S_n} 
\nc{\TmSm}{\T_m\setminus\S_m} 
\nc{\TXSX}{\T_X\setminus\S_X} 
\rnc{\S}{\mathcal S}

\nc{\T}{\mathcal T} 
\nc{\A}{\mathscr A} 
\nc{\B}{\mathcal B} 
\rnc{\P}{\mathcal P} 
\nc{\K}{\mathcal K}
\nc{\PB}{\mathcal{PB}} 
\nc{\rank}{\operatorname{rank}}

\nc{\mtt}{\!\!\!\mt\!\!\!}

\nc{\sub}{\subseteq}
\nc{\la}{\langle}
\nc{\ra}{\rangle}
\nc{\mt}{\mapsto}
\nc{\im}{\operatorname{im}}
\nc{\id}{\mathrm{id}}
\nc{\al}{\alpha}
\nc{\be}{\beta}
\nc{\ga}{\gamma}
\nc{\Ga}{\Gamma}
\nc{\de}{\delta}
\nc{\ka}{\kappa}
\nc{\lam}{\lambda}
\nc{\Lam}{\Lambda}
\nc{\si}{\sigma}
\nc{\Si}{\Sigma}
\nc{\oijn}{1\leq i<j\leq n}
\nc{\oijm}{1\leq i<j\leq m}

\nc{\comm}{\rightleftharpoons}
\nc{\AND}{\qquad\text{and}\qquad}

\nc{\bit}{\vspace{-3 truemm}\begin{itemize}}
\nc{\bitbmc}{\begin{itemize}\begin{multicols}}
%\nc{\bmc}{\vspace{-3 truemm}\begin{itemize}\begin{multicols}}
\nc{\bmc}{\begin{itemize}\begin{multicols}}
\nc{\emc}{\end{multicols}\end{itemize}\vspace{-3 truemm}}
\nc{\eit}{\end{itemize}\vspace{-3 truemm}}
\nc{\ben}{\vspace{-3 truemm}\begin{enumerate}}
\nc{\een}{\end{enumerate}\vspace{-3 truemm}}
\nc{\eitres}{\end{itemize}}

\nc{\set}[2]{\{ {#1} : {#2} \}} 
\nc{\bigset}[2]{\big\{ {#1}: {#2} \big\}} 
\nc{\Bigset}[2]{\left\{ \,{#1} :{#2}\, \right\}}

\nc{\pres}[2]{\la #1 : #2 \ra}
\nc{\bigpres}[2]{\big\la {#1} : {#2} \big\ra}
\nc{\Bigpres}[2]{\Big\la \,{#1}:{#2}\, \Big\ra}
\nc{\Biggpres}[2]{\Bigg\la {#1} : {#2} \Bigg\ra}
%\nc{\pres}[2]{\la {#1} \,|\, {#2} \ra}
%\nc{\bigpres}[2]{\big\la {#1} \,\big|\, {#2} \big\ra}
%\nc{\Bigpres}[2]{\Big\la \,{#1}\, \,\Big|\, \,{#2}\, \Big\ra}
%\nc{\Biggpres}[2]{\Bigg\la {#1} \,\Bigg|\, {#2} \Bigg\ra}

\nc{\pf}{\noindent{\bf Proof.}  }
\nc{\pfthm}[1]{\bigskip \noindent{\bf Proof of Theorem \ref{#1}.}  }
\nc{\epf}{\hfill$\Box$\bigskip}
\nc{\epfres}{\hfill$\Box$}
\nc{\pfnb}{\pf}
\nc{\epfnb}{\bigskip}
%\nc{\pfthm}[1]{\bigskip \noindent{\bf Proof of Theorem \ref{#1}}\,\,  } 
\nc{\pfprop}[1]{\bigskip \noindent{\bf Proof of Proposition \ref{#1}}\,\,  } 
%\nc{\pfthm}{\noindent{\bf Proof of Theorem \ref{mainthm} modulo Propositions \ref{prop1} and \ref{prop2}}\,\,  } 
\nc{\epfreseq}{\tag*{$\Box$}}

\nc{\uvert}[1]{\fill (#1,2)circle(.2);}
\rnc{\lvert}[1]{\fill (#1,0)circle(.2);}
\nc{\guvert}[1]{\fill[lightgray] (#1,2)circle(.2);}
\nc{\glvert}[1]{\fill[lightgray] (#1,0)circle(.2);}
\nc{\uvertx}[2]{\fill (#1,#2)circle(.2);}
\nc{\guvertx}[2]{\fill[lightgray] (#1,#2)circle(.2);}
\nc{\uvertxs}[2]{
\foreach \x in {#1}
{ \uvertx{\x}{#2}  }
}
\nc{\guvertxs}[2]{
\foreach \x in {#1}
{ \guvertx{\x}{#2}  }
}

\nc{\uvertth}[2]{\fill (#1,2)circle(#2);}
\nc{\lvertth}[2]{\fill (#1,0)circle(#2);}
\nc{\uvertths}[2]{
\foreach \x in {#1}
{ \uvertth{\x}{#2}  }
}
\nc{\lvertths}[2]{
\foreach \x in {#1}
{ \lvertth{\x}{#2}  }
}

\nc{\vertlabel}[2]{\draw(#1,2+.3)node{{\tiny $#2$}};}
\nc{\vertlabelh}[2]{\draw(#1,2+.4)node{{\tiny $#2$}};}
\nc{\vertlabelhh}[2]{\draw(#1,2+.6)node{{\tiny $#2$}};}
\nc{\vertlabelhhh}[2]{\draw(#1,2+.64)node{{\tiny $#2$}};}
\nc{\vertlabelup}[2]{\draw(#1,4+.6)node{{\tiny $#2$}};}
%\nc{\vertlabel}[2]{\draw(#1,2)node[above]{{\tiny $#2$}};}
\nc{\vertlabels}[1]{
{\foreach \x/\y in {#1}
{ \vertlabel{\x}{\y} }
}
}

\nc{\dvertlabel}[2]{\draw(#1,-.4)node{{\tiny $#2$}};}
\nc{\dvertlabels}[1]{
{\foreach \x/\y in {#1}
{ \dvertlabel{\x}{\y} }
}
}
\nc{\vertlabelsh}[1]{
{\foreach \x/\y in {#1}
{ \vertlabelh{\x}{\y} }
}
}
\nc{\vertlabelshh}[1]{
{\foreach \x/\y in {#1}
{ \vertlabelhh{\x}{\y} }
}
}
\nc{\vertlabelshhh}[1]{
{\foreach \x/\y in {#1}
{ \vertlabelhhh{\x}{\y} }
}
}

\nc{\vertlabelx}[3]{\draw(#1,2+#3+.6)node{{\tiny $#2$}};}
\nc{\vertlabelxs}[2]{
{\foreach \x/\y in {#1}
{ \vertlabelx{\x}{\y}{#2} }
}
}

\nc{\vertlabelupdash}[2]{\draw(#1,2.7)node{{\tiny $\phantom{'}#2'$}};}
\nc{\vertlabelupdashess}[1]{
{\foreach \x/\y in {#1}
{\vertlabelupdash{\x}{\y}}
}
}

\nc{\vertlabeldn}[2]{\draw(#1,0-.6)node{{\tiny $\phantom{'}#2'$}};}
\nc{\vertlabeldnph}[2]{\draw(#1,0-.6)node{{\tiny $\phantom{'#2'}$}};}

\nc{\vertlabelups}[1]{
{\foreach \x in {#1}
{\vertlabel{\x}{\x}}
}
}
\nc{\vertlabeldns}[1]{
{\foreach \x in {#1}
{\vertlabeldn{\x}{\x}}
}
}
\nc{\vertlabeldnsph}[1]{
{\foreach \x in {#1}
{\vertlabeldnph{\x}{\x}}
}
}

%\nc{\dotsup}[2]{\draw [dotted] (#1+.6,2)--(#2-.6,2);} \nc{\dotsdn}[2]{\draw [dotted] (#1+.6,0)--(#2-.6,0);}
\nc{\dotsup}[2]{\node()at(#1*.5+#2*.5,2){\tiny$\cdots$};} \nc{\dotsdn}[2]{\node()at(#1*.5+#2*.5,0){\tiny$\cdots$};}

\nc{\dotsupx}[3]{\draw [dotted] (#1+.6,#3)--(#2-.6,#3);}

\nc{\dotsups}[1]{\foreach \x/\y in {#1}
{ \dotsup{\x}{\y} }
}
\nc{\dotsupxs}[2]{\foreach \x/\y in {#1}
{ \dotsupx{\x}{\y}{#2} }
}
\nc{\dotsdns}[1]{\foreach \x/\y in {#1}
{ \dotsdn{\x}{\y} }
}

\nc{\nodropcustpartn}[3]{
\begin{tikzpicture}[scale=.3]
\foreach \x in {#1}
{ \uvert{\x}  }
\foreach \x in {#2}
{ \lvert{\x}  }
#3 \end{tikzpicture}
}

\nc{\custpartn}[3]{{\lower2 ex\hbox{
\begin{tikzpicture}[scale=.3]
\foreach \x in {#1}
{ \uvert{\x}  }
\foreach \x in {#2}
{ \lvert{\x}  }
#3 \end{tikzpicture}
}}}

\nc{\smcustpartn}[3]{{\lower0.7 ex\hbox{
\begin{tikzpicture}[scale=.2]
\foreach \x in {#1}
{ \uvert{\x}  }
\foreach \x in {#2}
{ \lvert{\x}  }
#3 \end{tikzpicture}
}}}

\nc{\dropcustpartn}[3]{{\lower5.2 ex\hbox{
\begin{tikzpicture}[scale=.3]
\foreach \x in {#1}
{ \uvert{\x}  }
\foreach \x in {#2}
{ \lvert{\x}  }
#3 \end{tikzpicture}
}}}

\nc{\dropcustpartnx}[4]{{\lower#4 ex\hbox{
\begin{tikzpicture}[scale=.4]
\foreach \x in {#1}
{ \uvert{\x}  }
\foreach \x in {#2}
{ \lvert{\x}  }
#3 \end{tikzpicture}
}}}

\nc{\dropcustpartnxy}[3]{\dropcustpartnx{#1}{#2}{#3}{4.6}}

\nc{\uvertsm}[1]{\fill (#1,2)circle(.15);}
\nc{\lvertsm}[1]{\fill (#1,0)circle(.15);}
\nc{\vertsm}[2]{\fill (#1,#2)circle(.15);}

\nc{\bigdropcustpartn}[3]{{\lower6.93 ex\hbox{
\begin{tikzpicture}[scale=.6]
\foreach \x in {#1}
{ \uvertsm{\x}  }
\foreach \x in {#2}
{ \lvertsm{\x}  }
#3 \end{tikzpicture}
}}}

\nc{\gcustpartn}[5]{{\lower1.4 ex\hbox{
\begin{tikzpicture}[scale=.3]
\foreach \x in {#1}
{ \uvert{\x}  }
\foreach \x in {#2}
{ \guvert{\x}  }
\foreach \x in {#3}
{ \lvert{\x}  }
\foreach \x in {#4}
{ \glvert{\x}  }
#5 \end{tikzpicture}
}}}

\nc{\gcustpartndash}[5]{{\lower6.97 ex\hbox{
\begin{tikzpicture}[scale=.3]
\foreach \x in {#1}
{ \uvert{\x}  }
\foreach \x in {#2}
{ \guvert{\x}  }
\foreach \x in {#3}
{ \lvert{\x}  }
\foreach \x in {#4}
{ \glvert{\x}  }
#5 \end{tikzpicture}
}}}

\nc{\stline}[2]{\draw(#1,2)--(#2,0);}
\nc{\stlines}[1]{
{\foreach \x/\y in {#1}
{ \stline{\x}{\y} }
}
}

\nc{\uarcs}[1]{
{\foreach \x/\y in {#1}
{ \uarc{\x}{\y} }
}
}

\nc{\darcs}[1]{
{\foreach \x/\y in {#1}
{ \darc{\x}{\y} }
}
}

\nc{\stlinests}[1]{
{\foreach \x/\y in {#1}
{ \stlinest{\x}{\y} }
}
}

\nc{\stlineds}[1]{
{\foreach \x/\y in {#1}
{ \stlined{\x}{\y} }
}
}

\nc{\gstline}[2]{\draw[lightgray](#1,2)--(#2,0);}
\nc{\gstlines}[1]{
{\foreach \x/\y in {#1}
{ \gstline{\x}{\y} }
}
}

\nc{\gstlinex}[3]{\draw[lightgray](#1,2+#3)--(#2,0+#3);}
\nc{\gstlinexs}[2]{
{\foreach \x/\y in {#1}
{ \gstlinex{\x}{\y}{#2} }
}
}

\nc{\stlinex}[3]{\draw(#1,2+#3)--(#2,0+#3);}
\nc{\stlinexs}[2]{
{\foreach \x/\y in {#1}
{ \stlinex{\x}{\y}{#2} }
}
}

\nc{\darcx}[3]{\draw(#1,0)arc(180:90:#3) (#1+#3,#3)--(#2-#3,#3) (#2-#3,#3) arc(90:0:#3);}
\nc{\darc}[2]{\darcx{#1}{#2}{.4}}
\nc{\uarcx}[3]{\draw(#1,2)arc(180:270:#3) (#1+#3,2-#3)--(#2-#3,2-#3) (#2-#3,2-#3) arc(270:360:#3);}
\nc{\uarc}[2]{\uarcx{#1}{#2}{.4}}

\nc{\darcxx}[4]{\draw(#1,0+#4)arc(180:90:#3) (#1+#3,#3+#4)--(#2-#3,#3+#4) (#2-#3,#3+#4) arc(90:0:#3);}
\nc{\uarcxx}[4]{\draw(#1,2+#4)arc(180:270:#3) (#1+#3,2-#3+#4)--(#2-#3,2-#3+#4) (#2-#3,2-#3+#4) arc(270:360:#3);}

\makeatletter
\newcommand\footnoteref[1]{\protected@xdef\@thefnmark{\ref{#1}}\@footnotemark}
\makeatother

\newcounter{theorem}
%\numberwithin{equation}{section}
\numberwithin{theorem}{section}

\newtheorem{thm}[theorem]{Theorem}
\newtheorem{lemma}[theorem]{Lemma}
\newtheorem{cor}[theorem]{Corollary}
\newtheorem{prop}[theorem]{Proposition}

\theoremstyle{definition}

\newtheorem{rem}[theorem]{Remark}
\newtheorem{defn}[theorem]{Definition}
\newtheorem{eg}[theorem]{Example}
\newtheorem{ass}[theorem]{Assumption}

\title{Presentations for (singular) partition monoids: \\ a new approach}

\date{}

\author{
James East\\
{\footnotesize \emph{Centre for Research in Mathematics; School of Computing, Engineering and Mathematics}}\\
{\footnotesize \emph{Western Sydney University, Locked Bag 1797, Penrith NSW 2751, Australia}}\\
{\footnotesize {\tt J.East\,@\,WesternSydney.edu.au}}
}

\maketitle

%\vspace{-0.5cm}

\begin{abstract}
We give new, short proofs of the presentations for the partition monoid and its singular ideal originally given in the author's 2011 papers in J.~Alg.~and I.J.A.C.

{\it Keywords}: Partition monoid, Singular ideal, Presentations.

MSC: 20M05; 20M20.
\end{abstract}

\section{Introduction}\label{sect:intro}

Partition monoids arise in the construction of partition algebras \cite{Martin1994,Jones1994_2} as twisted semigroup algebras \cite{Wilcox2007,HR2005}.  Although the partition algebras (and other related diagram algebras) were originally introduced in the context of theoretical physics and representation theory, their applications have been many and varied; see for example the recent surveys \cite{Koenig2008,Martin2008}.  In particular, diagram monoids have played an increasingly important role in semigroup theory in the last two decades; see for example \cite{EF,FL2011,LF2006,ACHLV2015,ADV2012_2,Auinger2014}, and especially \cite{EastGray} for an extensive list of references.  This semigroup approach was first employed by Wilcox in the context of cellular algebras \cite{Wilcox2007}, and was also instrumental in providing the first full proof \cite{JEgrpm} of a presentation (by generators and relations) for the partition algebras; the presentation was first stated in \cite{HR2005}.  The method in \cite{JEgrpm} was to first obtain a presentation for the partition monoid (stated in Theorem \ref{thm:Pn} below), and then apply a general result on twisted semigroup algebras (which was also proved in \cite{JEgrpm}).

Taking classical work of Howie on singular transformation semigroups \cite{Howie1966} as inspiration, the singular ideals of partition monoids were investigated in \cite{JEpnsn}, the main result being a presentation by generators and relations (stated in Theorem~\ref{thm:PnSn} below); the generating set consists of idempotents and is of minimal possible size.  (See also \cite{Maltcev2007}, where a similar presentation was given for the singular part of the Brauer monoid.)  The article \cite{JEpnsn} unlocked some intriguing combinatorial properties of the partition monoids, and paved the way for several further studies, most notably \cite{EastGray}, which concerns the question of (minimal) idempotent generation of arbitrary ideals in several natural families of diagram monoids.  The proofs of Theorems~\ref{thm:PnSn} and \ref{thm:Pn} given in \cite{JEgrpm,JEpnsn} relied on several previous results \cite{JEinsn,JEtnsn,Popova,FitzGerald2003} to obtain initial (highly complex) presentations,
%(involving five families of generators and fifteen families of relations), 
which were subsequently reduced (using lengthy sequences of Tietze transformations); the articles \cite{JEgrpm,JEpnsn} have a combined length of 58 pages.
% to an elegant presentation utilising a minimal-size all-idempotent generating set.  (See also \cite{Maltcev2007}, where a similar presentation was given for the singular part of the Brauer monoid.)

The purpose of the current article is to give shorter, more direct, and conceptually simpler proofs of Theorems \ref{thm:PnSn} and \ref{thm:Pn}.  
%
%Apart from quoting two results (from \cite{JEinsn2} and \cite{FitzGerald2003}), the current article is entirely self-contained: we give the required definitions and state the main result in Section~\ref{sect:prelim}, and then prove it in Section~\ref{sect:proof}.
%
Apart from quoting two results (Theorem 2.1 of \cite{JEinsn2} and Theorem 2 of \cite{FitzGerald2003}), the current article is entirely self-contained: we give the required definitions and state the main results in Section~\ref{sect:prelim}, and then prove them in Sections~\ref{sect:proof} and \ref{sect:Pn}.
%
%The current proof utilises \cite[Theorem 2.1]{JEinsn2} on singular symmetric inverse monoids, and \cite[Theorem 2]{FitzGerald2003} on idempotents in dual symmetric inverse monoids.  Apart from quoting these two results, the current article is entirely self-contained: we give the required definitions and state the main result in Section~\ref{sect:prelim}, and then prove it in Section~\ref{sect:proof}.
%
%We also hope the techniques introduced in this article will be of use to other investigators working on presentations.
%
%While the main results here are not new, we 
We believe that the new, efficient techniques for working with presentations are of independent interest, as may be the normal forms given in Proposition \ref{prop:tut3}.

\section{Preliminaries and statement of the main results}\label{sect:prelim}

Fix an integer $n\geq2$ (everything is trivial if $n\leq1$), and write $\bn=\{1,\ldots,n\}$ and $\bn'=\{1',\ldots,n'\}$.  The \emph{partition monoid of degree $n$}, denoted $\P_n$, consists of all set partitions of $\bn\cup\bn'$, under a product described below.
Such a partition $\al\in\P_n$ may be represented graphically.  We draw vertices $1,\ldots,n$ on an upper row (increasing from left to right) with $1',\ldots,n'$ directly below, and add edges so that connected components of the graph correspond to the blocks of $\al$.  For example, the partition
$\al=\big\{ \{1,4\},\{2,3,4',5'\},\{5,6\},\{1',3',6'\},\{2'\}\big\}\in\P_6$
is represented by the graph $\textpartnx6{\uarcx14{.6}\uarcx23{.3}\uarcx56{.3}\darc13\darcx36{.6}\darcx45{.3}\stline34}$.

The product of two partitions $\al,\be\in\P_n$ is calculated as follows.  We first stack (graphs representing)~$\al$ and $\be$ so that lower vertices $1',\ldots,n'$ of $\al$ are identified with upper vertices $1,\ldots,n$ of $\be$.  
%
%The product $\al\be$ is then (represented by) the graph on connected components of this graph are then constructed, and we finally delete the middle row; the resulting graph is the product $\al\be\in\P_n$.  
%
The connected components of this graph are then constructed, and we finally delete the middle row; the resulting graph is the product $\al\be\in\P_n$.  
Here is an example calculation with $\al,\be\in\P_6$:
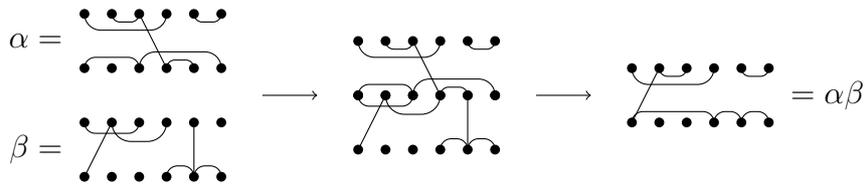
\[
\begin{tikzpicture}[scale=.36]
\begin{scope}[shift={(0,0)}]	
\buvs{1,...,6}
\blvs{1,...,6}
\uarcx14{.6}
\uarcx23{.3}
\uarcx56{.3}
\darc13
\darcx36{.6}
\darcx45{.3}
\stline34
\draw(0.6,1)node[left]{$\al=$};
\draw[->](7.5,-1)--(9.5,-1);
\end{scope}
\begin{scope}[shift={(0,-4)}]	
\buvs{1,...,6}
\blvs{1,...,6}
\uarc13
\uarcx24{.7}
\darc45
\darc56
\stline21
\stline55
\draw(0.6,1)node[left]{$\be=$};
\end{scope}
\begin{scope}[shift={(10,-1)}]	
\buvs{1,...,6}
\blvs{1,...,6}
\uarcx14{.6}
\uarcx23{.3}
\uarcx56{.3}
\darc13
\darcx36{.6}
\darcx45{.3}
\stline34
\draw[->](7.5,0)--(9.5,0);
\end{scope}
\begin{scope}[shift={(10,-3)}]	
\buvs{1,...,6}
\blvs{1,...,6}
\uarc13
\uarcx24{.7}
\darc45
\darc56
\stline21
\stline55
\end{scope}
\begin{scope}[shift={(20,-2)}]	
\buvs{1,...,6}
\blvs{1,...,6}
\uarcx14{.6}
\uarcx23{.3}
\uarcx56{.3}
\darc14
\darc45
\darc56
\stline21
\draw(6.4,1)node[right]{$=\al\be$};
\end{scope}
\end{tikzpicture}
%\vspace{-5mm}
\]
The operation is associative, so $\P_n$ is a semigroup: in fact, a monoid, with identity $1=\custpartn{1,2,4}{1,2,4}{\dotsups{2/4}\dotsdns{2/4}\stlines{1/1,2/2,4/4}}$.

A block of a partition $\al\in\P_n$ is called a \emph{transversal block} if it has non-trivial intersection with both~$\bn$ and $\bn'$; otherwise, it is a \emph{non-transversal block} (which may be an \emph{upper} or \emph{lower} non-transversal block, with obvious meanings).  The \emph{rank} of $\al$, denoted $\rank(\al)$, is the number of transversal blocks of $\al$.  For example, $\al\in\P_6$ as above has transversal block $\{2,3,4',5'\}$ (so $\rank(\al)=1$), upper non-transversal blocks $\{1,4\}$, $\{5,6\}$, and lower non-transversal blocks $\{1',3',6'\}$, $\{2'\}$.

For $\al\in\P_n$ and $x\in\bn\cup\bn'$, we write $[x]_\al$ for the block of $\al$ containing $x$.  We then define the \emph{domain} and \emph{codomain} of $\al$ to be the sets
\begin{align*}
\dom(\al) = \set{x\in\bn}{[x]_\al\cap\bn'\not=\emptyset} &\AND
\codom(\al) = \set{x\in\bn}{[x']_\al\cap\bn\not=\emptyset}.
\intertext{We also define the \emph{kernel} and \emph{cokernel} of $\al$ to be the equivalences}
\ker(\al) = \bigset{(x,y)\in\bn\times\bn}{[x]_\al=[y]_\al} &\AND
\coker(\al) = \bigset{(x,y)\in\bn\times\bn}{[x']_\al=[y']_\al}.
\end{align*}
For example, with $\al\in\P_6$ as as above, 
\begin{gather*}
\dom(\al)=\{2,3\} \COMMA \codom(\al) = \{4,5\} \COMMA
\ker(\al) = (1,4\mid2,3\mid5,6)  \COMMA \coker(\al) = (1,3,6\mid2\mid4,5),
\end{gather*}
using an obvious notation for equivalences.

It is immediate from the definitions that the following hold for all $\al,\be\in\P_n$:
\[
\begin{array}{rclcrcl}
%\rank(\al\be)\leq\min(\rank(\al),\rank(\be)) \COMMA
\dom(\al\be) \sub \dom(\al) \COMMA
\ker(\al\be)\sp \ker(\al) \COMMA
\codom(\al\be) \sub \codom(\be) \COMMA
\coker(\al\be)\sp \coker(\be).
\end{array}
\]
Also, $\rank(\al\be\ga)\leq\rank(\be)$ for all $\al,\be,\ga\in\P_n$.

If $\al\in\P_n$, we will write
\[
\al=\partnlist{A_1}{A_r}{C_1}{C_p}{B_1}{B_r}{D_1}{D_q}
\]
%to indicate that $\al$ has transversal blocks $A_i\cup B_i'$ (for $i\in\br$), upper non-transversal blocks $C_j$ (for~${j\in\bp}$), and lower non-transversal blocks $D_k'$ (for $k\in\bq$).  (For $A\sub\bn$, we write~${A'=\set{a'}{a\in A}}$.)  
to indicate that $\al$ has transversal blocks $A_i\cup B_i'$ (for $1\leq i\leq r$), upper non-transversal blocks~$C_j$ (for $1\leq j\leq p$), and lower non-transversal blocks $D_k'$ (for $1\leq k\leq q$).  (For $A\sub\bn$, we write~${A'=\set{a'}{a\in A}}$.)  
If $\al$ (as above) has no non-transversal blocks, we will write
$
{\al=\big({A_1\atop B_1} \divider {\cdots\atop\cdots} \divider {A_r\atop B_r}\big)}
$.
%\[
%\al=\partnlistnodefshort{A_1}{A_r}{B_1}{B_r}.
%\]
Such a partition is also known as a \emph{block bijection}, and the set $\J_n$ of all such block bijections is (isomorphic to) the \emph{dual symmetric inverse monoid of degree $n$}; see \cite{JEgrpm,FL1998}.  
If each block of $\al$ (as above) intersects~$\bn$ in at most one point and also intersects $\bn'$ in at most one point, we will write
$
{\al=\big[{a_1\atop b_1} \divider {\cdots\atop\cdots} \divider {a_r\atop b_r}\big]}
$
%\[
%\al=\partialpermshort{a_1}{a_r}{b_1}{b_r} %=\partnlist{a_1}{a_r}{c_1}{c_p}{B_1}{B_r}{D_1}{D_q}
%\]
to indicate that $\al$ has transversal blocks $\{a_i,b_i'\}$ (for $1\leq i\leq r$),  all other blocks being singletons.  Such a partition is also known as a \emph{partial permutation}, and the set $\I_n$ of all such partial permutations is (isomorphic to) the \emph{symmetric inverse monoid}; see \cite{JEpnsn,Lipscombe1996}.  As usual, with $\al\in\I_n$ as above, we will write $a_i\al=b_i$ for each $i$.  The group of units of $\P_n$ is (isomorphic to) the \emph{symmetric group} $\S_n=\J_n\cap\I_n$.  
A permutation $\pi$ of $\bn$ is identified with the partition $\big({1\atop 1\pi} \divider {\cdots\atop\cdots} \divider {n\atop n\pi}\big) = \big[{1\atop 1\pi} \divider {\cdots\atop\cdots} \divider {n\atop n\pi}\big]$.
%\[
%\partialpermshort1n{1\pi}{n\pi}=\partnlistnodefshort1n{1\pi}{n\pi}.
%\]
%
The set $\PnSn=\set{\al\in\P_n}{\rank(\al)<n}$ of non-invertible (i.e., \emph{singular}) partitions is a subsemigroup (indeed, an ideal) of $\P_n$.  
%
%A presentation was given for this ideal in \cite{JEpnsn}, and it is the goal of the current article to provide a much shorter (and more direct) proof of this presentation.  Our method relies on a recent presentation for the singular ideal $\InSn$ of the symmetric inverse monoid \cite{JEinsn2}.  

In order to state the above-mentioned presentations for $\PnSn$ and $\P_n$ from \cite{JEpnsn,JEgrpm}, we first fix the notation we will be using for presentations.  
%
%Let $X$ be an alphabet, and write $X^+$ for the free semigroup over $X$.  For a subset $R\sub X^+\times X^+$, we write $R^\sharp$ for the congruence on $X^+$ generated by $R$.  We say that a semigroup $S$ has \emph{(semigroup) presentation} $\pres XR$ if $S\cong X^+/R^\sharp$ or, equivalently, if there is an epimorphism $\phi:X^+\to S$ with $\ker\phi=R^\sharp$; in this case, we say $S$ has presentation $\pres XR$ \emph{via $\phi$}.  
%
Let $X$ be an alphabet, and denote by $X^+$ (resp., $X^*$) the free semigroup (resp., free monoid) on~$X$.  If~$R\sub X^+\times X^+$ (resp., $R\sub X^*\times X^*$), we denote by $R^\sharp$ the congruence on $X^+$ (resp., $X^*$) generated by $R$.  We say a semigroup (resp., monoid) $S$ has \emph{semigroup} (resp., \emph{monoid}) \emph{presentation} $\pres{X\!}{\!R}$ if~${S\cong X^+/R^\sharp}$ (resp., $S\cong X^*/R^\sharp$) or, equivalently, if there is an epimorphism ${X^{+}}\to S$ (resp., $X^*\to S$) with kernel $R^\sharp$.  If $\phi$ is such an epimorphism, we say $S$ has \emph{presentation $\pres{X\!}{\!R}$ via $\phi$}.  A relation $(w_1,w_2)\in R$ will usually be displayed as an equation: $w_1=w_2$.  We will always be careful to specify whether a given presentation is a semigroup or monoid presentation.  We denote the \emph{empty word} (over any alphabet) by $1$, so $X^*=X^+\cup\{1\}$.  If $A$ is a subset of a semigroup~$S$, then $\la A\ra$ always denotes the \emph{subsemigroup} generated by $A$.

%All presentations considered here are semigroup presentations.  Generation is always understood to be within the variety of semigroups; if $A$ is a subset of a semigroup $S$.
%A word $x_i\cdots x_j$ is understood 

%\section{Statement of the main result}\label{sect:statement}

For $1\leq r\leq n$ and $\oijn$, define partitions
\[
\eb_r = \custpartn{1,3,4,5,7}{1,3,4,5,7}{\dotsups{1/3,5/7}\dotsdns{1/3,5/7}\stlines{1/1,3/3,5/5,7/7}\vertlabelshh{1/1,4/r,7/n}} \AND
 \tb_{ij} = \tb_{ji} = \custpartn{1,3,4,5,7,8,9,11}{1,3,4,5,7,8,9,11}{\dotsups{1/3,5/7,9/11}\dotsdns{1/3,5/7,9/11}\stlines{1/1,3/3,4/4,5/5,7/7,8/8,9/9,11/11}\uarc48 \darc48\vertlabelshh{1/1,4/i,8/j,11/n}}.
\]
%The reason for the (over-line) notation will become clear shortly.
%
%For $\ve\in\Eq_n$ with equivalence classes $A_1,\ldots,A_r$, we define
%\[
%\tb_\ve = \partnlistnodefshort{A_1}{A_r}{A_1}{A_r}.
%\]
%For 
%\[
%\al=\partnlist{A_1}{A_r}{C_1}{C_p}{B_1}{B_r}{D_1}{D_q}\in \P_n,
%\]
%define
%\[
%\widehat\al = \partnlistnodefshort{a_1}{a_r}{b_1}{b_r}\in\I_n,
%\]
%where $a_i=\min(A_i)$ and $b_i=\min(B_i)$ for each $i$.
%
%\ms
%\begin{prop}\label{prop:normalform}
%Let $\al\in\P_n$ and write $\ve=\ker(\al)$ and $\eta=\coker(\al)$.  Then $\al = \tb_\ve \widehat\al \tb_\eta$.  Further, if $\be\in\P_n$, then $\al=\be$ if and only if $\ker(\be)=\ve$, $\coker(\be)=\eta$ and $\widehat\be=\widehat\al$.
%\end{prop}
%
Consider alphabets $E=\set{e_r}{1\leq r\leq n}$ and $T=\set{t_{ij}}{\oijn}$, and define a (semigroup) homomorphism
\[
{\phi:(E\cup T)^+\to\PnSn}
\]
by $e_r\phi =  \eb_r$ and $t_{ij}\phi =  \tb_{ij}$ for each $1\leq r\leq n$ and $\oijn$.
We will use symmetric notation when referring to the letters from $T$, so we write $t_{ij}=t_{ji}$ for all $\oijn$.
Consider the relations
\begin{align}
\tag{R1} e_i^2 &= e_i  &&\text{{for all $i$}}\\
\tag{R2} e_ie_j &= e_je_i   &&\text{{for distinct $i,j$}}\\
\tag{R3} t_{ij}^2 &= t_{ij}   &&\text{{for all $i,j$}}\\
\tag{R4} t_{ij}t_{kl} &= t_{kl}t_{ij}  &&\text{{for all $i,j,k,l$}}\\
\tag{R5} t_{ij}t_{jk} &= t_{jk}t_{ki}    &&\text{{for distinct $i,j,k$}}\\
\tag{R6} t_{ij}e_k &= e_kt_{ij}  &&\text{{if $k\not\in\{i,j\}$}}\\
\tag{R7} t_{ij}e_kt_{ij} &=  t_{ij} &&\text{{if $k\in\{i,j\}$}}\\
\tag{R8} e_kt_{ij}e_k &=  e_k &&\text{{if $k\in\{i,j\}$}}\\
\tag{R9} e_kt_{ki}e_it_{ij}e_jt_{jk}e_k &= e_kt_{kj}e_jt_{ji}e_it_{ik}e_k  &&\text{{for distinct $i,j,k$}}\\% &&\text{\emph{for all\, $i,j,k$}}\\
\tag{R10} e_kt_{ki}e_it_{ij}e_jt_{jl}e_lt_{lk}e_k &= e_kt_{kl}e_lt_{li}e_it_{ij}e_jt_{jk}e_k &&\text{{for distinct $i,j,k,l$.}}
\end{align}
Here is the first of our main results; it originally appeared in \cite[Theorem 46]{JEpnsn}.
%As stated above, one of our two main goals is to provide a short proof of the following result, which was originally proved in \cite[Theorem 46]{JEpnsn}.

\ms
\begin{thm}\label{thm:PnSn}
The semigroup $\PnSn$ has semigroup presentation $\pres{E\cup T}{\text{\emph{(R1--R10)}}}$ via $\phi$. 
\end{thm}

%Write
%\[
%\E = E\phi = \set{\ve_i}{1\leq i\leq n} \AND \T = T\phi = \set{\tau_{ij}}{\oijn}.
%\]
%
%
%
%\newpage

In order to state the second main result, define partitions
\[
\eb=\eb_1= \custpartn{1,2,3,5}{1,2,3,5}{\dotsups{3/5}\dotsdns{3/5}\stlines{2/2,3/3,5/5}\vertlabelshh{1/1,5/n}}\COMMA 
\tb=\tb_1= \custpartn{1,2,3,5}{1,2,3,5}{\dotsups{3/5}\dotsdns{3/5}\stlines{1/1,2/2,3/3,5/5}\vertlabelshh{1/1,5/n}\arcup12\arcdn12}\COMMA
\sb_i\Phi= \custpartn{1,3,4,5,6,8}{1,3,4,5,6,8}{\dotsups{1/3,6/8}\dotsdns{1/3,6/8}\stlines{1/1,3/3,4/5,5/4,6/6,8/8}\vertlabelshh{1/1,4/i,8/n}} \qquad\text{for $1\leq i\leq n-1$.}
\]
Consider the alphabet $S=\set{s_i}{1\leq i\leq n-1}$, and define a (monoid) homomorphism
\[
\Phi:(S\cup\{e,t\})^*\to\P_n
\]
by $e\Phi=\eb$, $t\Phi=\tb$ and $s_i\Phi=\sb_i$ for each $i$.  Consider the relations
\begin{align}
\tag{R11} s_i^2 &= 1 &&\hspace{-4 cm}\text{for all $i$}\\
\tag{R12} s_is_j &= s_js_i &&\hspace{-4 cm}\text{if $|i-j|>1$}\\
\tag{R13} s_is_js_i &= s_js_is_j &&\hspace{-4 cm}\text{if $|i-j|=1$}\\
\tag{R14} e^2 = e &= e te\\
\tag{R15} t^2=t=te t&=ts_1=s_1t\\
\tag{R16} e s_i &= s_i e &&\hspace{-4 cm}\text{if $i\geq2$}\\
\tag{R17} t s_i &= s_i t &&\hspace{-4 cm}\text{if $i\geq3$}\\
\tag{R18} s_1e s_1e = e &s_1e s_1 = e s_1e\\
\tag{R19} t s_2 t s_2 &= s_2 t s_2 t \\
\tag{R20} t (s_2s_3s_1s_2) t (s_2s_3s_1s_2) &= (s_2s_3s_1s_2) t (s_2s_3s_1s_2) t\\
\tag{R21} t(s_2s_1e s_1s_2) &= (s_2s_1es_1s_2)t. 
\end{align}
Here is our second main result; this originally appeared in \cite[Theorem 32]{JEgrpm}.

\ms
\begin{thm}\label{thm:Pn}
The monoid $\P_n$ has monoid presentation $\pres{S\cup\{e,t\}}{\text{\emph{(R11--R21)}}}$ via $\Phi$.
\end{thm}

We conclude this section by stating the presentation
%A crucial ingredient in our proof of Theorem \ref{thm:PnSn} is a presentation 
for $\InSn$ from \cite{JEinsn2}.  (Recall that $\I_n$ was defined above.)
With this in mind, for $i,j\in\bn$ with $i\not=j$, define $\fb_{ij}=\eb_i\tb_{ij}\eb_j$  (using symmetric notation for $\tb_{ij}=\tb_{ji}$).  Note that $\fb_{ij}\not=\fb_{ji}$; rather,
\[
\fb_{ij} = \begin{cases}
\custpartn{1,3,4,5,7,8,9,11}{1,3,4,5,7,8,9,11}{\dotsups{1/3,5/7,9/11}\dotsdns{1/3,5/7,9/11}\stlines{1/1,3/3,5/5,7/7,8/4,9/9,11/11}\vertlabelshh{1/1,4/i,8/j,11/n}}
&\text{if $i<j$}
\\~\\
\custpartn{1,3,4,5,7,8,9,11}{1,3,4,5,7,8,9,11}{\dotsups{1/3,5/7,9/11}\dotsdns{1/3,5/7,9/11}\stlines{1/1,3/3,5/5,7/7,4/8,9/9,11/11}\vertlabelshh{1/1,4/j,8/i,11/n}}
&\text{if $j<i$.}
\end{cases}
\]
%
%For $\oijn$, define 
%\[
%\fb_{ij} = \custpartn{1,3,4,5,7,8,9,11}{1,3,4,5,7,8,9,11}{\dotsups{1/3,5/7,9/11}\dotsdns{1/3,5/7,9/11}\stlines{1/1,3/3,5/5,7/7,8/4,9/9,11/11}\vertlabelshh{1/1,4/i,8/j,11/n}}
%\AND
%\fb_{ji} = \custpartn{1,3,4,5,7,8,9,11}{1,3,4,5,7,8,9,11}{\dotsups{1/3,5/7,9/11}\dotsdns{1/3,5/7,9/11}\stlines{1/1,3/3,5/5,7/7,4/8,9/9,11/11}\vertlabelshh{1/1,4/i,8/j,11/n}}.
%\]
%%
%%$\fb_{ij}\in\I_n$ by
%%\[
%%k\fb_{ij} = \begin{cases}
%%k &\text{if $k\in\{i,j\}^c$}\\
%%i &\text{if $k=j$}\\
%%- &\text{if $k=i$.}
%%\end{cases}
%%\]
%%The reason for the notation will become clear later.  
%It is easy to check diagrammatically that $\fb_{ij}=\eb_i\tb_{ij}\eb_j$ for all $i,j\in\bn$ with $i\not=j$ (using symmetric notation for $\tb_{ij}=\tb_{ji}$).

Define an alphabet $F=\set{f_{ij}}{i,j\in\bn,\ i\not=j}$, and consider the relations
\begin{align}
\label{F1}\tag{F1} f_{ij}f_{ji}f_{ij} &= f_{ij} &&\text{{for distinct $i,j$}}\\ 
\label{F2}\tag{F2} f_{ij}^3=f_{ij}^2 &= f_{ji}^2 &&\text{{for distinct $i,j$}}\\ 
\label{F3}\tag{F3} f_{ij}f_{kl} &= f_{kl}f_{ij} &&\text{{for distinct $i,j,k,l$}}\\ 
\label{F4}\tag{F4} f_{ij}f_{ji} &= f_{ik}f_{ki} &&\text{{for distinct $i,j,k$}}\\ 
\label{F5}\tag{F5} f_{ij}f_{ik}=f_{jk}f_{ij} &= f_{ik}f_{jk} &&\text{{for distinct $i,j,k$}}\\ 
\label{F6}\tag{F6} f_{ki}f_{ij}f_{jk} &= f_{kj}f_{ji}f_{ik} &&\text{{for distinct $i,j,k$}}\\ 
\label{F7}\tag{F7} f_{ki}f_{ij}f_{jk}f_{kl} &= f_{kl}f_{li}f_{ij}f_{jl} &&\text{{for distinct $i,j,k,l$.}}
\end{align}
The next result is \cite[Theorem 2.1]{JEinsn2}.

\ms
\begin{thm}\label{thm:InSn}
The semigroup $\InSn$ has semigroup presentation $\pres{F}{\text{\emph{(F1--F7)}}}$ via $f_{ij}\mt\fb_{ij}$. \epfres
\end{thm}

\section{Presentation for $\PnSn$}\label{sect:proof}

%This section is devoted to a proof of Theorem \ref{thm:PnSn}.  
There are two components to the proof of Theorem \ref{thm:PnSn}: that the map $\phi:(E\cup T)^+\to\PnSn$ is an epimorphism; and that $\ker\phi$ is generated by relations (R1--R10).  Proposition~\ref{prop:surjective} accomplishes the first task, and the remainder of the section is devoted to the second.  
%For the proof of Proposition~\ref{prop:surjective}, and for later use, we first describe the join semilattice of equivalence relations.

We write $\Eq_n$ for the set of all equivalence relations on $\bn$, which we regard as a semilattice (monoid of commuting idempotents)
%(monoid of commuting idempotents) 
under $\vee$; the \emph{join}, $\ve\vee\eta$, of two equivalences $\ve,\eta\in\Eq_n$ is defined to be the smallest equivalence containing $\ve\cup\eta$ (see \cite{FL1998,FitzGerald2003}).

For a subset $A\sub\bn$ with $\bn\sm A=\{i_1,\ldots,i_k\}$, we write 
%\[
%\tb_A = \left( \begin{array}{c|c|c|c} 
%\!\! A & i_1 & \cdots & i_k\ \ \\ 
%\!\! A & i_1 & \cdots & i_k \ \
%\end{array} \!\!\! \right).
%\]
$
\tb_A = \big(  {A\atop A} \divider {i_1\atop i_1} \divider {\cdots\atop \cdots} \divider {i_k\atop i_k} \big)
$.
Note that $\tb_A=1$ if $|A|\leq1$.  
%For $\oijn$, we write $\tb_{ij}=\tb_{ji}=\tb_{\{i,j\}}$.  
Note that if $i,j\in\bn$ with $i\not=j$, then $\tb_{\{i,j\}}=\tb_{ij}$ in the notation of the previous section.
Note that if $A=\{a_1,a_2,\ldots,a_r\}$ with $r=|A|\geq2$, then $\tb_A=\tb_{a_1a_2}\tb_{a_2a_3}\cdots \tb_{a_{r-1}a_r}=\tb_{a_1a_2}\tb_{a_1a_3}\cdots \tb_{a_1a_r}$.  
For an equivalence $\ve\in\Eq_n$ with equivalence classes $A_1,\ldots,A_k$, we write
%\[
%\tb_\ve=\partnlistnodefshort{A_1}{A_k}{A_1}{A_k}=\tb_{A_1}\cdots\tb_{A_k}.
%\]
$
\tb_\ve=\big( {A_1\atop A_1} \divider {\cdots\atop\cdots} \divider {A_k\atop A_k} \big)=\tb_{A_1}\cdots\tb_{A_k}
$.
In fact, the set $\set{\tb_\ve}{\ve\in\Eq_n}$ is equal to $E(\J_n)$, the semilattice of idempotents of (the isomorphic copy of) the dual symmetric inverse monoid $\J_n\sub\P_n$; see \cite{FL1998}.
%; the \emph{join} $\ve\vee\eta$ of two equivalences $\ve,\eta\in\Eq_n$ is defined to be the least equivalence containing $\ve\cup\eta$ (see \cite{FL1998,FitzGerald2003}).  
In particular, $E(\J_n)$ is isomorphic to the semilattice $(\Eq_n,\vee)$; so $\tb_\ve\tb_\eta=\tb_{\ve\vee\eta}$ for all $\ve,\eta\in\Eq_n$.  By the discussion above, we see that $E(\J_n)$ is generated (as a monoid) by the set ${\set{\tb_{ij}}{\oijn}}$.  Theorem 2 of \cite{FitzGerald2003} says that $E(\J_n)$ has (monoid) presentation $\pres T{\text{(R3--R5)}}$ via the map $t_{ij}\mt\tb_{ij}$.

\ms
\begin{prop}\label{prop:surjective}
The semigroup $\PnSn$ is generated by the set $\set{\eb_r}{1\leq r\leq n}\cup\set{\tb_{ij}}{\oijn}$.  In particular, the map $\phi:(E\cup T)^+\to\PnSn$ is surjective.
\end{prop}

\pf Let $\al\in\PnSn$, and write
\[
\al=\partnlist{A_1}{A_r}{C_1}{C_p}{B_1}{B_r}{D_1}{D_q},
\]
noting that $r=\rank(\al)\leq n-1$.  For each $1\leq i\leq r$, choose some $a_i\in A_i$ and $b_i\in B_i$.  Then clearly $\al=\be\ga\de$, where
$\be=\big( {A_1\atop A_1} \divider {\cdots\atop\cdots} \divider {A_r\atop A_r} \divider {C_1\atop C_1} \divider {\cdots\atop\cdots} \divider {C_p\atop C_p} \big)$, $\ga=\big[{a_1\atop b_1} \divider {\cdots\atop\cdots} \divider {a_r\atop b_r}\big]$ and $\de=\big( {B_1\atop B_1} \divider {\cdots\atop\cdots} \divider {B_r\atop B_r} \divider {D_1\atop D_1} \divider {\cdots\atop\cdots} \divider {D_q\atop D_q} \big)$.
%\[
%\be=\partnlistnodef{A_1}{A_r}{C_1}{C_p}{A_1}{A_r}{C_1}{C_p} \COMMA
%\ga=\partialpermshort{a_1}{a_r}{b_1}{b_r} \COMMA
%\de=\partnlistnodef{B_1}{B_r}{D_1}{D_q}{B_1}{B_r}{D_1}{D_q}.
%\]
Note that $\be=\tb_{\ker(\al)}$ and $\de=\tb_{\coker(\al)}$, so that $\be,\de\in\pres{\tb_{ij}}{\oijn}\cup\{1\}$, by the discussion before the statement of the proposition.  By Theorem \ref{thm:InSn}, $\ga\in\InSn$ is a (non-empty) product of terms of the form $\fb_{ij}=\eb_i\tb_{ij}\eb_j$. \epf

\ms
\begin{rem}
Proposition \ref{prop:surjective} was also proven in \cite[Theorem 9]{EF}, using a classical result of Howie on transformation semigroups \cite{Howie1966}.  Since $\set{\eb_r}{1\leq r\leq n}\sub E(\I_n)$ and $\set{\tb_{ij}}{\oijn}\sub E(\J_n)$, it follows that $\la E(\P_n)\ra=\la E(\I_n)\cup E(\J_n)\ra=\{1\}\cup(\PnSn)$.
\end{rem}

Now that we know $\phi$ is surjective, it remains to show that $\ker\phi={\sim}$, where $\sim$ is the congruence on $(E\cup T)^+$ generated by the relations (R1--R10).  For $w\in(E\cup T)^+$, we write $\wb=w\phi$.  Even though the empty word $1$ does not belong to $(E\cup T)^+$, we will also write $\overline1=1$ and $1\sim1$.

\ms
\begin{lemma}
We have ${\sim}\sub\ker\phi$.
\end{lemma}

\pf This follows by a simple diagrammatic check that $\phi$ preserves the relations (R1--R10).  We do this for (R7) in Figure \ref{fig:R7}, and leave the rest for the reader. \epf 

\begin{figure}[h]
\begin{center}
\begin{tikzpicture}[scale=.41]
\begin{scope}[shift={(0,2)}]	
\uverts{1,3,4,5,7,8,9,11}
\lverts{1,3,4,5,7,8,9,11}
\stlines{1/1,3/3,4/4,5/5,7/7,8/8,9/9,11/11}
\uarc48 
\darc48
\dotsups{1/3,5/7,9/11}
\dotsdns{1/3,5/7,9/11}
\vertlabelshh{1/1,4/i,8/j,11/n}
\end{scope}
\begin{scope}[shift={(0,0)}]	
\uverts{1,3,4,5,7,8,9,11}
\lverts{1,3,4,5,7,8,9,11}
\stlines{1/1,3/3,5/5,7/7,8/8,9/9,11/11}
\dotsups{1/3,5/7,9/11}
\dotsdns{1/3,5/7,9/11}
\draw(14,1)node{$=$};
\end{scope}
\begin{scope}[shift={(0,-2)}]	
\uverts{1,3,4,5,7,8,9,11}
\lverts{1,3,4,5,7,8,9,11}
\stlines{1/1,3/3,4/4,5/5,7/7,8/8,9/9,11/11}
\uarc48 
\darc48
\dotsups{1/3,5/7,9/11}
\dotsdns{1/3,5/7,9/11}
\end{scope}
\begin{scope}[shift={(16,2)}]	
\uverts{1,3,4,5,7,8,9,11}
\lverts{1,3,4,5,7,8,9,11}
\stlines{1/1,3/3,4/4,5/5,7/7,8/8,9/9,11/11}
\uarc48 
\darc48
\dotsups{1/3,5/7,9/11}
\dotsdns{1/3,5/7,9/11}
\vertlabelshh{1/1,4/i,8/j,11/n}
\end{scope}
\begin{scope}[shift={(16,0)}]	
\uverts{1,3,4,5,7,8,9,11}
\lverts{1,3,4,5,7,8,9,11}
\stlines{1/1,3/3,5/5,7/7,4/4,9/9,11/11}
\dotsups{1/3,5/7,9/11}
\dotsdns{1/3,5/7,9/11}
\draw(14,1)node{$=$};
\end{scope}
\begin{scope}[shift={(16,-2)}]	
\uverts{1,3,4,5,7,8,9,11}
\lverts{1,3,4,5,7,8,9,11}
\stlines{1/1,3/3,4/4,5/5,7/7,8/8,9/9,11/11}
\uarc48 
\darc48
\dotsups{1/3,5/7,9/11}
\dotsdns{1/3,5/7,9/11}
\end{scope}
\begin{scope}[shift={(32,0)}]	
\uverts{1,3,4,5,7,8,9,11}
\lverts{1,3,4,5,7,8,9,11}
\stlines{1/1,3/3,4/4,5/5,7/7,8/8,9/9,11/11}
\uarc48 
\darc48
\dotsups{1/3,5/7,9/11}
\dotsdns{1/3,5/7,9/11}
\vertlabelshh{1/1,4/i,8/j,11/n}
\end{scope}
\end{tikzpicture}
\end{center}
\vspace{-5mm}
\caption{Diagrammatic proof of relation (R7): $\tb_{ij}\eb_i\tb_{ij}=\tb_{ij}\eb_j\tb_{ij}=\tb_{ij}$.}
\label{fig:R7}
\end{figure}

For $i,j\in\bn$ with $i\not=j$, define the word $z_{ij}=e_it_{ij}e_j\in(E\cup T)^+$, and write
$
Z = \set{z_{ij}}{i,j\in\bn,\ i\not=j}
$
for the set of all such words.  Note that $\zb_{ij}=z_{ij}\phi=\fb_{ij}\in\I_n$.
% (the $\fb_{ij}$ we defined after the statement of Theorem \ref{thm:PnSn}).
%
The words $z_{ij}$ will play a crucial role in what follows.  We first prove a number of basic relations satisfied among them.
%Write $\zeta_{ij}=z_{ij}\phi$ for each $i,j$ (draw).  So
%\[
%\InSn=\la\Z\ra \qquad\text{where}\qquad \Z=Z\phi=\set{\zeta_{ij}}{i,j\in\bn,\ i\not=j}.
%\]
%{\red $Z$ or $F$?}

\ms
\begin{lemma}\label{lem:EZ}
Let $i,j,k\in\bn$ be distinct.  Then
\[
\text{\emph{
(i) $e_iz_{ij} \sim z_{ij} \sim z_{ij}e_j$; \ \ \ \  
(ii) $e_jz_{ij}\sim z_{ij}e_i\sim e_ie_j \sim z_{ij}^2\sim z_{ji}^2$; \ \ \ \ 
(iii) $z_{ij}z_{ji}\sim e_i$; \ \ \ \ 
%(iv) $z_{ij}^2\sim e_ie_j$; \ \ \ \ 
(iv) $e_kz_{ij}\sim z_{ij}e_k$.}}
\]
\end{lemma}

\pf Part (i) follows immediately from (R1).  For (ii), note that $e_jz_{ij} = e_j e_it_{ij}e_j \sim e_ie_jt_{ij}e_j \sim e_ie_j$, by (R2) and (R8).  A similar calculation gives $z_{ij}e_i\sim e_ie_j$.  Together with (R8), it also follows that $z_{ij}^2=z_{ij}e_i t_{ij}e_j \sim e_ie_j t_{ij}e_j\sim e_ie_j$; similarly, $z_{ji}^2\sim e_ie_j$, completing the proof of (ii).  For (iii), we have $z_{ij}z_{ji} = e_it_{ij}e_je_jt_{ij}e_i \sim e_it_{ij}e_jt_{ij}e_i \sim e_it_{ij}e_i \sim e_i$, by (R1), (R7) and (R8).  Finally, (iv) follows immediately from (R2) and (R6). \epf

Recall that for $\al\in\I_n$ and $x\in\dom(\al)$, we write $x\al$ for the (unique) element of $\codom(\al)$ such that $\{x,(x\al)'\}$ is a block of $\al$.  Recall also that we are using symmetric notation for the $t_{ij}=t_{ji}$.

\ms
\begin{lemma}\label{lem:tijzkl}
Let $i,j,k,l\in\bn$ with $i\not=j$, $k\not=l$ and $k\not\in\{i,j\}$.  Then
\[
t_{ij} z_{kl} \sim z_{kl}t_{i\zb_{kl},j\zb_{kl}} =
\begin{cases}
z_{kl}t_{ij} &\text{if $l\not\in\{i,j\}$}\\
z_{kl}t_{kj} &\text{if $l=i$}\\
z_{kl}t_{ik} &\text{if $l=j$.}
\end{cases}
\]
\end{lemma}

\pf The case in which $l\not\in\{i,j\}$ follows immediately from (R4) and (R6).  If $l=i$, then
\[
t_{ij} z_{kl} = t_{lj}e_kt_{kl}e_l \sim e_kt_{lj}t_{kl}e_l =  e_kt_{jl}t_{lk}e_l \sim e_kt_{lk}t_{kj}e_l \sim e_kt_{kl}e_lt_{kj} = z_{kl}t_{kj},
\]
by (R5) and (R6).  The $l=j$ case is similar. \epf

%Note that Lemma \ref{lem:tijzkl} says that $t_{ij} z_{kl} \sim z_{kl}t_{i\zb_{kl},j\zb_{kl}} $ for all $i,j,k,l\in\bn$ with the stated constraints.

Recall that $\la Z\ra$ denotes the subsemigroup of $(E\cup T)^+$ generated by $Z$.

\ms
\begin{cor}\label{cor:twwt}
If $w\in\la Z\ra$ and $i,j\in\dom(\wb)$ with $i\not=j$, then $t_{ij}w\sim wt_{i\wb,j\wb}$.
\end{cor}

\pf This follows immediately from Lemma \ref{lem:tijzkl}, after writing $w=z_{k_1l_1}\cdots z_{k_sl_s}$. \epf

The next result shows that (modulo the relations (R1--R10)) the words $z_{ij}$ satisfy the defining relations for $\InSn$ from Theorem \ref{thm:InSn}.

\ms
\begin{lemma}\label{lem:InSn}
For distinct $i,j,k,l\in\bn$, we have
\begin{align}
\label{Z1}\tag{Z1} z_{ij}z_{ji}z_{ij} &\sim z_{ij} \\ 
\label{Z2}\tag{Z2} z_{ij}^3\sim z_{ij}^2 &\sim z_{ji}^2 \\ 
\label{Z3}\tag{Z3} z_{ij}z_{kl} &\sim z_{kl}z_{ij} \\ 
\label{Z4}\tag{Z4} z_{ij}z_{ji} &\sim z_{ik}z_{ki} \\ 
\label{Z5}\tag{Z5} z_{ij}z_{ik}\sim z_{jk}z_{ij} &\sim z_{ik}z_{jk}\\ 
\label{Z6}\tag{Z6} z_{ki}z_{ij}z_{jk} &\sim z_{kj}z_{ji}z_{ik} \\ 
\label{Z7}\tag{Z7} z_{ki}z_{ij}z_{jk}z_{kl} &\sim z_{kl}z_{li}z_{ij}z_{jl} . 
\end{align}
\end{lemma}

\pf For (Z1), we have $z_{ij}z_{ji}z_{ij}\sim e_iz_{ij}\sim z_{ij}$, by Lemma \ref{lem:EZ}(iii) and (i).  For (Z2), note that $z_{ij}^2\sim z_{ji}^2\sim e_ie_j$ by Lemma \ref{lem:EZ}(ii); it also follows that $z_{ij}^3\sim e_ie_jz_{ij}\sim e_ie_ie_j \sim e_ie_j$ by Lemma \ref{lem:EZ}(ii) and (R1), completing the proof of (Z2).  Relation (Z3) follows immediately from (R2), (R4) and~(R6).  Relation (Z4) follows immediately from Lemma \ref{lem:EZ}(iii).  
For (Z5), first note that $z_{ij}z_{ik} \sim z_{ij}e_iz_{ik} \sim e_ie_jz_{ik}\sim e_je_iz_{ik}\sim e_jz_{ik}$, by Lemma \ref{lem:EZ}(i) and (ii), and (R2).  We also have
\begin{align*}
z_{jk}z_{ij} = e_jt_{jk}e_k e_it_{ij}e_j &\sim e_ie_jt_{jk} t_{ij}e_je_k &&\hspace{-6cm}\text{by (R2) and (R6)}\\
&=  e_ie_jt_{kj} t_{ji}e_je_k \\
&\sim e_ie_jt_{ji} t_{ik}e_je_k &&\hspace{-6cm}\text{by (R5)}\\
&\sim e_ie_jt_{ji} e_jt_{ik}e_k &&\hspace{-6cm}\text{by (R6)}\\
&\sim e_ie_jt_{ik}e_k &&\hspace{-6cm}\text{by (R8)}\\
&\sim e_je_it_{ik}e_k = e_jz_{ik} &&\hspace{-6cm}\text{by (R2).}
%\end{align*}
\intertext{By Lemma \ref{lem:EZ}(i), (ii) and (iv), %together with (R2), we have
$z_{ik}z_{jk} \sim z_{ik}e_kz_{jk} \sim z_{ik}e_je_k \sim e_jz_{ik}e_k \sim \sim e_jz_{ik}$, completing the proof of (Z5).  For (Z6), we have}
%\begin{align*}
%z_{ki}z_{ij}z_{jk} = e_kt_{ki}e_i  e_it_{ij}e_j  e_jt_{jk}e_k 
%&\sim e_kt_{ki}e_i t_{ij}e_j t_{jk}e_k &&\text{by (R1)}\\
%&\sim e_kt_{kj}e_j t_{ji}e_i t_{ik}e_k &&\text{by (R9)}\\
%&\sim e_kt_{kj}e_j  e_jt_{ji}e_i  e_it_{ik}e_k = z_{kj}z_{ji}z_{ik} &&\text{by (R1).}
%\end{align*}
%\[
z_{ki}z_{ij}z_{jk} = e_kt_{ki}e_i  e_it_{ij}e_j  e_jt_{jk}e_k 
\sim e_k&t_{ki}e_i t_{ij}e_j t_{jk}e_k 
\sim e_kt_{kj}e_j t_{ji}e_i t_{ik}e_k 
\sim e_kt_{kj}e_j  e_jt_{ji}e_i  e_it_{ik}e_k = z_{kj}z_{ji}z_{ik} ,
%\]
\intertext{by (R1) and (R9).
Finally, for (Z7), first observe that by a similar calculation to that just carried out, (R10) gives $z_{ki}z_{ij}z_{jl}z_{lk} \sim z_{kl}z_{li}z_{ij}z_{jk}$.
We then have}
%\begin{align*}
z_{ki}z_{ij}z_{jk}z_{kl} 
&\sim z_{ki}z_{ij}z_{ji}z_{ij}z_{jk}z_{kl} &&\hspace{-6cm}\text{by (Z1)}\\
%&\sim z_{ki}z_{ij}z_{jk}z_{kj}z_{jk}z_{kl} &&\text{by (Z4)}\\
&\sim z_{ki}z_{ij}z_{jl}z_{lj}z_{jk}z_{kl} &&\hspace{-6cm}\text{by (Z4)}\\
&\sim z_{ki}z_{ij}z_{jl}z_{lk}z_{kj}z_{jl} &&\hspace{-6cm}\text{by (Z6)}\\
&\sim z_{kl}z_{li}z_{ij}z_{jk}z_{kj}z_{jl} &&\hspace{-6cm}\text{by the observation}\\
&\sim z_{kl}z_{li}z_{ij}z_{jl}z_{lj}z_{jl} &&\hspace{-6cm}\text{by (Z4)}\\
&\sim z_{kl}z_{li}z_{ij}z_{jl} &&\hspace{-6cm}\text{by (Z1).}
\end{align*}
This completes the proof. \epf

\ms
\begin{cor}\label{cor:InSn}
If $u,v\in\la Z\ra$, then $\ub=\vb \implies u\sim v$.
\end{cor}

\pf Write $u=z_{i_1j_1}\cdots z_{i_sj_s}$ and $v=z_{k_1l_1}\cdots z_{k_tl_t}$, and suppose $\ub=\vb$.  Then
\[
\fb_{i_1j_1}\cdots \fb_{i_sj_s} = \zb_{i_1j_1}\cdots \zb_{i_sj_s} = \ub=\vb = \zb_{k_1l_1}\cdots \zb_{k_tl_t} = \fb_{k_1l_1}\cdots \fb_{k_tl_t}.
\]
By Theorem \ref{thm:InSn}, the word $f_{i_1j_1}\cdots f_{i_sj_s}$ may be transformed into $f_{k_1l_1}\cdots f_{k_tl_t}$ using relations (F1--F7).  By Lemma \ref{lem:InSn}, this transformation leads to a transformation of $u$ into $v$ using (Z1--Z7).  \epf

%The next result follows from the simple observation that (R1--R10) contains defining relations in a presentation \cite{FitzGerald2003} for $\Eq_n$.

The next result follows from \cite[Theorem 2]{FitzGerald2003}, which (as noted above) states that $E(\J_n)\cong(\Eq_n,\vee)$ has monoid presentation $\pres T{\text{(R3--R5)}}$ via $t_{ij}\mt\tb_{ij}$.

\ms
\begin{lemma}\label{lem:T}
If $u,v\in T^*$, then $\ub=\vb \implies u\sim v$. \epfres
\end{lemma}

%\pf This follows from \cite[Theorem 2]{FitzGerald2003}, which (as noted above) states that $E(\J_n)\cong(\Eq_n,\vee)$ has monoid presentation $\pres T{\text{R3--R5)}}$ via $t_{ij}\mt\tb_{ij}$. \epf

%\pf 
%%{\bf (Z1):} 
%For (Z1), we have $z_{ij}z_{ji}z_{ij} = e_it_{ij}e_je_jt_{ij}e_ie_it_{ij}e_j \sim e_it_{ij}e_jt_{ij}e_it_{ij}e_j \sim e_it_{ij}e_jt_{ij}e_j \sim e_it_{ij}e_j = z_{ij}$, by (R1) and (R7).
%
%For (Z2), first note that $z_{ij}^2 = e_it_{ij}e_j\cdot e_it_{ij}e_j \sim e_it_{ij}e_i\cdot e_jt_{ij}e_j \sim e_ie_j$, by (R2) and (R8).  Together with (R2), it follows, $z_{ji}^2\sim e_je_i\sim e_ie_j\sim z_{ij}^2$.  We also have $z_{ij}^3\sim e_ie_j\cdot$

%\ms
%\begin{lemma}\label{lem:zijtkl}
%Let $i,j\in\bn$ with $i\not=j$, and let $k,l\in\{j\}^c$ with $k\not=l$.  Then
%\[
%z_{ij} t_{kl} \sim \begin{cases}
%t_{kl}z_{ij} &\text{if $i\not\in\{k,l\}$}\\
%t_{kj}z_{ij} &\text{if $i=l$}\\
%t_{jl}z_{ij} &\text{if $i=k$.}
%\end{cases}
%\]
%\end{lemma}
%
%\pf The case in which $i\not\in\{k,l\}$ follows immediately from (R4) and (R6).  If $i=l$, then
%\[
%z_{ij}t_{kl} = z_{ij}t_{ki} = e_it_{ij}e_j t_{ki} \sim e_it_{ij}t_{ki}e_j = e_it_{ji}t_{ik}e_j \sim e_it_{kj}t_{ji}e_j  \sim t_{kj}e_it_{ij}e_j = t_{kj}z_{ij},
%\]
%by (R5) and (R6).  The $i=k$ case is similar. \epf

In order to complete the proof of Theorem \ref{thm:PnSn}, we aim to show that any word over $E\cup T$ may be rewritten (using the relations) to take on a very specific form; see Proposition \ref{prop:tut3}, the proof of which requires the next three intermediate lemmas.

\ms
\begin{lemma}\label{lem:w123}
If $w\in(E\cup T)^+$, then $w\sim w_1w_2w_3$ for some $w_1,w_3\in T^*$ and $w_2\in\la Z\ra$.
\end{lemma}

\pf We prove this by induction on $\ell(w)$, the length of the word $w$.  Suppose first that $\ell(w)=1$.  If $w=e_i$ for some $i$, then $w\sim z_{ij}z_{ji}$ for any $j\in\bn\sm\{i\}$, by Lemma \ref{lem:EZ}(iii), and we are done (with $w_1=w_3=1$ and $w_2=z_{ij}z_{ji}$).
%\[
%w=e_i \sim e_it_{ij}e_i \sim e_it_{ij}e_jt_{ij}e_i \sim e_it_{ij}e_je_jt_{ij}e_i = z_{ij}z_{ji},
%\]
%by (R8), (R7) and (R1), respectively.  
If $w=t_{ij}$ for some $\oijn$, then
$w=t_{ij} \sim t_{ij}e_it_{ij} \sim t_{ij}z_{ij}z_{ji}t_{ij}$,
by (R7) and Lemma \ref{lem:EZ}(iii), and we are done (with $w_1=w_3=t_{ij}$ and $w_2=z_{ij}z_{ji}$).  

Now suppose $\ell(w)\geq2$, and write $w=ux$, where $u\in(E\cup T)^+$ and $x\in E\cup T$.  By an inductive hypothesis, $u\sim u_1u_2u_3$ for some $u_1,u_3\in T^*$ and $u_2\in\la Z\ra$.  If $x\in T$, then $w\sim u_1u_2u_3x$, and we are done (with $w_1=u_1$, $w_2=u_2$ and $w_3=u_3x$).  So suppose $x=e_i\in E$, and write $u_3=t_{k_1l_1}\cdots t_{k_sl_s}$.  

{\bf Case 1.}  If $i\not\in\{k_1,\ldots,k_s,l_1,\ldots,l_s\}$, then $u_3e_i\sim e_iu_3\sim z_{ij}z_{ji}u_3$ for any $j\in\bn\sm\{i\}$, by (R6) and Lemma~\ref{lem:EZ}(iii), so $w\sim u_1u_2z_{ij}z_{ji}u_3$, and we are done (with $w_1=u_1$, $w_2=u_2z_{ij}z_{ji}$ and $w_3=u_3$).

{\bf Case 2.}  Now suppose $i\in\{k_1,\ldots,k_s,l_1,\ldots,l_s\}$.  By (R3), (R4), and the symmetrical notation for the $t_{kl}=t_{lk}$, we may assume that in fact $u_3=t_{il_1}\cdots t_{il_r}\cdot t_{k_{r+1}l_{r+1}}\cdots t_{k_sl_s}$ with $r\geq1$ and $l_1<\cdots<l_r$.  Now $t_{il_1}\cdots t_{il_r} \sim t_{il_1} t_{l_1l_2}\cdots t_{l_{r-1}l_r}$, by Lemma \ref{lem:T}.  It follows that 
\[
w \sim u_1u_2  t_{il_1}  t_{l_1l_2}\cdots t_{l_{r-1}l_r} t_{k_{r+1}l_{r+1}}\cdots t_{k_sl_s} e_i \sim u_1\cdot u_2  t_{il_1}  e_i \cdot u_4,
\]
where $u_4= t_{l_1l_2}\cdots t_{l_{r-1}l_r} t_{k_{r+1}l_{r+1}}\cdots t_{k_sl_s}$.  For simplicity, we will write $j=l_1$.  Since $u_1,u_4\in T^*$, the proof will be complete if we can show that $u_2t_{ij}e_i\sim v_1v_2v_3$ for some $v_1,v_3\in T^*$ and $v_2\in\la Z\ra$.  

{\bf Subcase 2.1.}  Suppose first that $i,j\in\codom(\ub_2)$, and let $k,l\in\dom(\ub_2)$ be such that $k\ub_2=i$ and $l\ub_2=j$.
%put $k=i\ub_2^{-1}$ and $l=j\ub_2^{-1}$, 
Then $t_{kl}u_2\sim u_2t_{ij}$, by Corollary \ref{cor:twwt}, whence $u_2t_{ij}e_i \sim t_{kl}u_2e_i \sim t_{kl}u_2z_{ij}z_{ji}$, and we are done (with $v_1= t_{kl}$, $v_2=u_2z_{ij}z_{ji}$ and $v_3=1$).

{\bf Subcase 2.2.}  Next suppose $i\not\in\codom(\ub_2)$.  Then $u_2\sim u_2z_{ij}z_{ji}\sim u_2e_i$, by Corollary \ref{cor:InSn} and Lemma~\ref{lem:EZ}(iii).  But then, together with (R8), it follows that $u_2t_{ij}e_i \sim u_2e_it_{ij}e_i \sim u_2e_i\sim u_2$, and we are done (with $v_1=v_3=1$ and $v_2=u_2$).

%Then $u_2\sim u_2z_{ij}z_{ji}$, by Lemma \ref{lem:InSn}, so
%\[
%u_2t_{ij}e_i \sim u_2z_{ij}z_{ji}t_{ij}e_i = u_2 \cdot e_it_{ij}e_j\cdot e_jt_{ij}e_i \cdot t_{ij}\cdot e_i \sim u_2 \cdot e_it_{ij}e_j\cdot e_jt_{ij}e_i =u_2z_{ij}z_{ji} \sim u_2,
%\]
%by (R8), and we are done (with $v_1=v_3=1$ and $v_2=u_2$).

%{\bf Subcase 2.2.}  Next, suppose $i,j\in\codom(\ub_2)^c$.  Then $u_2\sim u_2z_{ji}^2$, so
%\[
%u_2t_{ij}e_i \sim u_2z_{ji}^2t_{ij}e_i = u_2 \cdot e_jt_{ij}e_i\cdot e_jt_{ij}e_i \cdot t_{ij}\cdot e_i \sim u_2 \cdot e_jt_{ij}e_i\cdot e_jt_{ij}e_i =u_2z_{ji}^2 \sim u_2,
%\]
%by (R8), and we are done (with $v_1=v_3=1$ and $v_2=u_2$).

{\bf Subcase 2.3.}  Finally, suppose $j\not\in\codom(\ub_2)$.  Then $u_2\sim u_2e_j$, as in the previous case, giving $u_2t_{ij}e_i \sim u_2e_jt_{ij}e_i = u_2z_{ji}$, and we are done (with $v_1=v_3=1$ and $v_2=u_2z_{ij}$). \epf

%\[
%u_2\sim u_2z_{ji}z_{ij} = u_2\cdot e_jt_{ij}e_i \cdot e_it_{ij}e_j \sim u_2\cdot e_jt_{ij}e_i t_{ij}e_j \sim u_2\cdot e_jt_{ij}e_j \sim u_2\cdot e_j,
%\]
%by Lemma \ref{lem:InSn}, (R1), (R7) and (R8).  It then follows that $u_2t_{ij}e_i \sim u_2e_jt_{ij}e_i = u_2z_{ji}$, and again we are done (with $v_1=v_3=1$ and $v_2=u_2z_{ji}$). \epf

To improve Lemma \ref{lem:w123}, we first define some words over $T$.
Consider a subset $A\sub\bn$.  If $|A|\leq 1$, then put $t_A=1$.  Otherwise, write $A=\{i_1,\ldots,i_k\}$ with $i_1<\cdots <i_k$, and define $t_A=t_{i_1i_2}t_{i_2i_3}\cdots t_{i_{k-1}i_k}$.  For an equivalence $\ve\in\Eq_n$ with equivalence classes $A_1,\ldots,A_r$ with $\min(A_1)<\cdots<\min(A_r)$, define $t_\ve=t_{A_1}\cdots t_{A_r}$.  Note that if $w\in T^*$ is such that
%\[
%\wb = \tb_\ve = \partnlistnodefshort{A_1}{A_r}{A_1}{A_r},
%\]
$
\wb = \tb_\ve = \big({A_1\atop A_1} \divider {\cdots\atop\cdots} \divider {A_r\atop A_r}\big)
$,
then $w\sim t_\ve$, by Lemma \ref{lem:T}.
%using relations (R3--R5), by \cite{FitzGerald2003}.  
For the proof of the next result, if $\ve\in\Eq_n$, we write $\bn/\ve$ for the set of all $\ve$-classes.

\ms
\begin{lemma}\label{lem:tut1}
Let $w\in(E\cup T)^+$, and put $\ve=\ker(\wb)$ and $\eta=\coker(\wb)$.  Then $w\sim t_\ve ut_\eta$ for some $u\in\la Z\ra$.
\end{lemma}

\pf By Lemmas \ref{lem:T} and \ref{lem:w123}, the set $\bigset{(\lam,u,\rho)\in\Eq_n\times\la Z\ra\times\Eq_n}{w\sim t_\lam ut_\rho}$ is non-empty.  Choose an element $(\lam,u,\rho)$ from this set such that $k=|\bn/\lam|+|\bn/\rho|$ is minimal.  
Suppose the $\lam$-classes and $\rho$-classes are $A_1,\ldots,A_p$ and $B_1,\ldots,B_q$, respectively (so $p=|\bn/\lam|$ and $q=|\bn/\rho|$).  Note that $\lam = \ker(\tb_\lam) \sub \ker(\tb_\lam \ub\tb_\rho) = \ker(\wb)=\ve$ and, similarly, $\rho\sub\eta$.  In particular, each $\ve$-class is a union of (one or more) $\lam$-classes, with a similar statement holding for $\eta$-classes and $\rho$-classes.  Let the $\ve$-classes and $\eta$-classes be $C_1,\ldots,C_s$ and $D_1,\ldots,D_t$, respectively (so $s=|\bn/\ve|$ and $t=|\bn/\eta|$), noting that $s\leq p$ and $t\leq q$.  In particular, $k=p+q\geq s+t$.  If $k=s+t$, then $p=s$ and $q=t$, so that $\lam=\ve$ and $\rho=\eta$, and the proof would be complete.  So suppose instead that $k>s+t$.  Note that one of the following statements must be true:
\bit
\item[(i)] there exist $a,b\in\dom(\ub)$ such that $(a,b)\in\lam$ and $(a\ub,b\ub)\not\in\rho$; or
\item[(ii)] there exist $a,b\in\dom(\ub)$ such that $(a,b)\not\in\lam$ and $(a\ub,b\ub)\in\rho$.
%\item $a,b\in\dom(\al)$ and $c,d\in\codom(\al)$ such that $a,b$ belong to the same $\lam$-class, $c,d$ belong to different $\rho$-classes, and $a\al=c$, $b\al=d$.
\eit
Indeed, if (i) and (ii) were both false, then we would have $\lam=\ker(\tb_\lam \ub\tb_\rho)=\ker(\wb)=\ve$ and, similarly, $\rho=\eta$, contradicting the assumption that $k>s+t$.  Suppose (i) is true (the other case being similar), and write $c=a\ub$ and $d=b\ub$.  Relabelling if necessary, we may assume that $a,b\in A_1$, $c\in B_1$ and $d\in B_2$.  Since $(a,b)\in\lam$, Lemma \ref{lem:T} gives $t_\lam\sim t_\lam t_{ab}$.  Let $\ve_{cd}\in\Eq_n$ be the equivalence whose only non-singleton equivalence class is $\{c,d\}$, and let $\ka=\ve_{cd}\vee\rho$.  Then Lemma \ref{lem:T} also gives $t_{cd}t_\rho=t_{\ve_{cd}}t_\rho\sim t_\ka$.  Together with Corollary \ref{cor:twwt}, it then follows that
$w \sim t_\lam t_{ab} u t_\rho \sim t_\lam u t_{cd} t_\rho \sim t_\lam u t_\ka$.
But $\ka$ has $q-1$ equivalence classes (i.e., $B_1\cup B_2,B_3,\ldots,B_q$), so $|\bn/\lam|+|\bn/\ka|=k-1$, contradicting the minimality of $k$. \epf

\ms
\begin{rem}
The proof of Lemma \ref{lem:tut1} is set out as a \emph{reductio} merely for convenience.  It is easily adapted to become constructive.  (The same is true of Lemma \ref{lem:tut2} and Proposition \ref{prop:tut3}).
\end{rem}

\ms
\begin{lemma}\label{lem:tut2}
Let $w\in(E\cup T)^+$, and put $\ve=\ker(\wb)$ and $\eta=\coker(\wb)$.  Then $w\sim t_\ve ut_\eta$ for some $u\in\la Z\ra$ with $\rank(\ub)=\rank(\wb)$.
%such that $\dom(\ub)$ (resp., $\codom(\ub)$) intersects each $\ve$-class (resp., $\eta$-class) in at most one point.
\end{lemma}

\pf Put $r=\rank(\wb)$.  By Lemma \ref{lem:tut1}, the set $\set{u\in\la Z\ra}{w\sim t_\ve ut_\eta}$ is non-empty.  Choose an element $u$ from this set with $\rank(\ub)$ minimal.  If $\rank(\ub)=r$, then we are done, so suppose otherwise.  Since $r=\rank(\wb)=\rank(\tb_\ve \ub\tb_\eta)\leq\rank(\ub)$, it follows that $\rank(\ub)>r$.
%
%$w\sim t_\ve ut_\eta$ for some $u\in\la Z\ra$.  Put $r=\rank(\wb)$ and $s=\rank(\ub)$, and suppose that $s=\min\set{\rank(\vb)}{v\in\la Z\ra,\ w\sim t_\ve vt_\eta}$. 
%
%note that $r=\rank(\overline{t_\ve vt_\eta})\leq\rank(\vb)$.  We proceed by induction on $\rank(\vb)$.  If $\rank(\vb)=r$, then we are already done, so suppose $\rank(\vb)>r$.
%If $\rank(\vb)=\rank(\wb)$, then we are done, so suppose $\rank(\vb)>\rank(\wb)$.  Put $r=\rank(\wb)$, and 
Suppose the $\ve$-classes and $\eta$-classes are $A_1,\ldots,A_p$ and $B_1,\ldots,B_q$, respectively, and 
%also that we have labelled these classes so 
that the transversal blocks of~$\wb$ are $A_1\cup B_1',\ldots,A_r\cup B_r'$.  
Note that $\dom(\ub)\sub A_1\cup\cdots\cup A_r$ and $\codom(\ub)\sub B_1\cup\cdots\cup B_r$.  
Since $\rank(\ub)>r$, we may assume (relabelling if necessary) that $|A_1\cap\dom(\ub)|\geq2$.  Let $i,j\in A_1\cap\dom(\ub)$ with $i\not=j$, and put $k=i\ub$ and $l=j\ub$, noting that $k,l\in B_1$.  It follows from Lemma \ref{lem:T} and Corollary~\ref{cor:twwt} that $t_\ve\sim t_\ve t_{ij}$, $t_\eta\sim t_{kl}t_\eta$, and $t_{ij}u\sim ut_{kl}$.  Together with (R7) and Lemma~\ref{lem:EZ}(iii), we then obtain
$
w\sim t_\ve ut_\eta \sim t_\ve t_{ij}ut_\eta \sim t_\ve t_{ij}e_it_{ij}ut_\eta \sim t_\ve t_{ij}e_iut_{kl}t_\eta \sim t_\ve(z_{ij}z_{ji}u)t_\eta.
$
But $\dom(\zb_{ij}\zb_{ji}\ub)=\dom(\ub)\sm\{i\}$, so that $\rank(\zb_{ij}\zb_{ji}\ub)=\rank(\ub)-1$, contradicting the minimality of $\rank(\ub)$. \epf

The next result gives a set of \emph{normal forms} for words over $E\cup T$, and is the final ingredient in the proof of Theorem \ref{thm:PnSn}.  For each $\al\in\InSn$, fix some $z_\al\in\la Z\ra$ such that $\zb_\al=z_\al\phi=\al$.  

\ms
\begin{prop}\label{prop:tut3}
Let $w\in(E\cup T)^+$, and write
\[
\wb = \partnlist{A_1}{A_r}{C_1}{C_p}{B_1}{B_r}{D_1}{D_q} \COMMA \ve=\ker(\wb) \COMMA \eta=\coker(\wb) \COMMA\al = \partialpermshort{a_1}{a_r}{b_1}{b_r},
\]
where $a_i=\min(A_i)$ and $b_i=\min(B_i)$ for each $i$.  Then $w\sim t_\ve z_\al t_\eta$.
%Then $w\sim t_\ve ut_\eta$ for some $u\in\la Z\ra$ with $\ub=\al$.
\end{prop}

\pf By Lemma \ref{lem:tut2}, the set $\set{u\in\la Z\ra}{\rank(\ub)=r,\ w\sim t_\ve ut_\eta}$ is non-empty.  Choose some $u$ from this set with $k=|{\dom(\ub)\cap\dom(\al)}|+|{\codom(\ub)\cap\codom(\al)}|$ maximal.  We claim that $k=2r$.  Indeed, suppose to the contrary that $k<2r$, and write $\ub=\big[{c_1\atop d_1} \divider {\cdots\atop\cdots} \divider {c_r\atop d_r}\big]$, where $c_i\in A_i$ and $d_i\in B_i$ for each $i$.  Since $k<2r$, it follows that $c_i\not=a_i$ or $d_i\not=b_i$ for some $i$.  We assume the former is the case (the latter is treated in similar fashion).  Relabelling if necessary, we may assume that $i=1$; for simplicity, we will write $a=a_1$ and $c=c_1$.  Since $a\not\in\dom(\ub)$, we have $u\sim z_{ac}z_{ca}u=e_at_{ac}e_cz_{ca}u\sim e_at_{ac}z_{ca}u$, by Corollary \ref{cor:InSn} and Lemma \ref{lem:EZ}(i).  Since $a,c\in A_1$, Lemma~\ref{lem:T} gives $t_\ve\sim t_\ve t_{ac}$.  Together with (R7), it follows that $w\sim t_\ve ut_\eta \sim t_\ve t_{ac} (e_at_{ac}z_{ca}u)t_\eta \sim t_\ve t_{ac} z_{ca}ut_\eta \sim t_\ve (z_{ca}u)t_\eta$.  Note that $\zb_{ca}\ub = \big[{a_1\atop d_1} \divider {c_2\atop d_2} \divider {\cdots\atop\cdots} \divider {c_r\atop d_r}\big]$.  But $|{\dom(\zb_{ca}\ub)\cap\dom(\al)}|+|{\codom(\zb_{ca}\ub)\cap\codom(\al)}|=k+1$, contradicting the maximality of~$k$.  This completes the proof of the claim.  It follows that $\ub=\al$, so Corollary \ref{cor:InSn} gives $u\sim z_\al$.  \epf

We now have all we need to complete the proof of the first main result.

\pfthm{thm:PnSn} It remains to check that $\ker\phi\sub{\sim}$, so suppose $w,v\in(E\cup T)^+$ are such that $\wb=\vb$.  Then $w\sim t_\ve z_\al t_\eta$, in the notation of Proposition \ref{prop:tut3}.  Since $\vb=\wb$, we also have $v\sim t_\ve z_\al t_\eta$, so that $w\sim v$. \epf

\section{Presentation for $\P_n$}\label{sect:Pn}

Again, the proof of Theorem \ref{thm:Pn} involves two steps: showing that the map $\Phi:(S\cup\{e,t\})^*\to\P_n$ is an epimorphism; and showing that $\ker\Phi$ is generated by relations (R11--R21).

Write $\approx$ for the congruence on $(S\cup\{e,t\})^*$ generated by (R11--R21).  Without causing confusion, we will write $\wb=w\Phi$ for any $w\in(S\cup\{e,t\})^*$.  For $w=s_{i_1}\cdots s_{i_k}\in S^*$, we will write $w^{-1}=s_{i_k}\cdots s_{i_1}$, noting that $ww^{-1}\approx w^{-1}w\approx 1$, by (R11).  Note also that the symmetric group $\S_n\sub\P_n$ has monoid presentation $\pres{S}{\text{(R11--R13)}}$ via $s_i\mt\sb_i$; see \cite{Moo}.

%\ms
%\begin{prop}
%The monoid $\P_n$ is generated (as a semigroup) by the set $\set{\sb_r}{1\leq r\leq n-1}\cup\{\eb,\tb\}$.  In particular, the map $\Phi:(S\cup\{e,t\})^*\to\P_n$ is surjective.
%\end{prop}

%\pf One may easily check, diagrammatically 

For $1\leq r\leq n$ and $\oijn$, define the words
\[
c_r=s_1\cdots s_{r-1} \COMMA \ep_r=c_i^{-1}ec_i \COMMA \tau_{ij} =\tau_{ji} = c_i^{-1}c_j^{-1}tc_jc_i.
\]
One may easily check diagrammatically that $\overline{\ep}_r=\ep_r\Phi=\eb_r$ and $\overline{\tau}_{ij}=\tau_{ij}\Phi=\tb_{ij}$.  In particular, $\im\Phi$ contains $\PnSn$, by Proposition \ref{prop:surjective}.  Since also $S\Phi=\set{\sb_r}{1\leq r\leq n-1}$ generates $\S_n$, it follows that $\Phi$ is surjective.  It is easy to check that each of relations (R11--R21) are preserved by $\Phi$.  So, to prove Theorem \ref{thm:Pn}, it remains to check that $\ker\Phi\sub{\approx}$, and the rest of this section is devoted to that task.  We begin with some simple properties of the words $\ep_r$ and $\tau_{ij}$.

\ms
\begin{lemma}\label{lem:ses}
If $1\leq r\leq n$ and $1\leq k\leq n-1$, then $s_k\ep_rs_k \approx \ep_{r\sb_k}$.
%\[
%s_k\ep_rs_k \approx \ep_{r\sb_k}=
%\begin{cases}
%\ep_{r-1} &\text{if $k=r-1$}\\
%\ep_{r+1} &\text{if $k=r$}\\
%\ep_r &\text{otherwise.}
%\end{cases}
%\]
\end{lemma}

\pf We must show that
\[
s_k\ep_rs_k \approx 
\begin{cases}
\ep_{r-1} &\text{if $k=r-1$}\\
\ep_{r+1} &\text{if $k=r$}\\
\ep_r &\text{otherwise.}
\end{cases}
\]
These follows quickly from (R11), (R16), and the easily checked facts that
\[
c_rs_k \approx 
\begin{cases}
s_{k+1}c_r &\text{if $k\leq r-2$}\\
c_{r-1} &\text{if $k=r-1$}\\
c_{r+1} &\text{if $k=r$}\\
s_kc_r &\text{if $k\geq r+1$}
\end{cases}
\AND
s_kc_r^{-1} \approx 
\begin{cases}
c_r^{-1} s_{k+1}&\text{if $k\leq r-2$}\\
c_{r-1}^{-1} &\text{if $k=r-1$}\\
c_{r+1}^{-1} &\text{if $k=r$}\\
c_r^{-1} s_k&\text{if $k\geq r+1$.}
\end{cases}
\]
For example, if $k=r-1$, then $s_k\ep_rs_k=s_{r-1}c_r^{-1}ec_rs_{r-1}\approx c_{r-1}^{-1}ec_{r-1}=\ep_{r-1}$, while if 
$k\leq r-2$, then
$s_k\ep_rs_k=s_kc_r^{-1}ec_rs_k \approx c_r^{-1}s_{k+1}es_{k+1}c_r \approx c_r^{-1}s_{k+1}s_{k+1}ec_r \approx c_r^{-1}ec_r=\ep_r$.~\epf
%\[
%\epfreseq
%s_k\ep_rs_k=s_kc_r^{-1}ec_rs_k \approx c_r^{-1}s_{k+1}es_{k+1}c_r \approx c_r^{-1}s_{k+1}s_{k+1}ec_r \approx c_r^{-1}ec_r=\ep_r.
%\]

\ms
\begin{lemma}[{cf.~\cite[p322]{FitzGerald2003}}]\label{lem:sts}
If $\oijn$ and $1\leq k\leq n-1$, then $s_k\tau_{ij}s_k \approx \tau_{i\sb_k,j\sb_k}$.
%\[
%s_k\tau_{ij}s_k \approx \tau_{i\sb_k,j\sb_k}=
%\begin{cases}
%\tau_{i-1,j} &\text{if $k=i-1$}\\
%\tau_{i+1,j} &\text{if $k=i<j-1$}\\
%\tau_{i,j-1} &\text{if $k=j-1>i$}\\
%\tau_{i,j+1} &\text{if $k=j$}\\
%\tau_{ij} &\text{otherwise.}
%\end{cases}
%\]
\end{lemma}

\pf This follows by a similar proof to that of Lemma \ref{lem:ses}, using the same relations satisfied by the $c_r,s_k$.  For example, if $k=j$, then
\[
\epfreseq
s_k\tau_{ij}s_k = s_j c_i^{-1}c_j^{-1}tc_jc_i s_j \approx  c_i^{-1}s_jc_j^{-1}tc_js_jc_i \approx c_i^{-1}c_{j+1}^{-1}tc_{j+1}c_i = \tau_{i,j+1}.
\]
%if $k=i<j-1$, then $s_k\tau_{ij}s_k = s_i c_i^{-1}c_j^{-1}tc_jc_i s_i = s_{i+1}^{-1}c_j^{-1}tc_jc_i s_{i+1}$. \epf

\ms
\begin{cor}\label{cor:wewwtw}
If $1\leq r\leq n$, $\oijn$, and $w\in S^*$, then $w^{-1}\ep_rw\approx\ep_{r\wb}$ and $w^{-1}\tau_{ij}w\approx\tau_{i\wb,j\wb}$.
\end{cor}

\pf This follows from Lemmas \ref{lem:ses} and \ref{lem:sts} and a simple induction on the length of $w$. \epf

We now aim to link the presentations $\pres{E\cup T}{\text{(R1--R10)}}$ and $\pres{S\cup\{e,t\}}{\text{(R11--R21)}}$ in a certain sense.  With this in mind, we define a map $\psi:(E\cup T)^+\to(S\cup\{e,t\})^*$ by $e_r\psi=\ep_r$ and $t_{ij}\psi=\tau_{ij}$.  Note that $\ub=\overline{u\psi}$ (i.e., $u\phi=u\psi\Phi$) for all $u\in(E\cup T)^+$.  Recall that $\sim$ is the congruence on $(E\cup T)^+$ generated by relations (R1--R10).

\ms
\begin{lemma}\label{lem:reduction1}
If $u,v\in(E\cup T)^+$, then $u\sim v \implies u\psi\approx v\psi$.
\end{lemma}

\pf We just need to check this for each relation $u=v$ from (R1--R10).  Relations (R1) and~(R3) follow immediately from (R11), (R14), (R15).  For (R2), first note that $\ep_1\ep_2=es_1es_1\approx s_1es_1e=\ep_2\ep_1$, by (R18).  Together with Corollary \ref{cor:wewwtw}, it follows that for any $w\in S^*$ with $1\wb=i$ and $2\wb=j$, $\ep_i\ep_j\approx w^{-1}\ep_1ww^{-1}\ep_2w \approx w^{-1}\ep_1\ep_2w \approx w^{-1}\ep_2\ep_1w \approx w^{-1}\ep_2ww^{-1}\ep_1w \approx \ep_j\ep_i$.  The other relations are all treated in similar fashion; in each case, we use Corollary \ref{cor:wewwtw} to reduce the calculation to a fixed set of values of the subscripts.  For example, for~(R9), taking $(i,j,k)=(1,2,3)$,
\begin{align*}
\ep_3 \tau_{31} \ep_1 \tau_{12} \ep_2 \tau_{23} \ep_3 &= (s_2s_1es_1s_2)(s_2s_1ts_1s_2)et(s_1es_1)(s_1s_2s_1ts_1s_2s_1)(s_2s_1es_1s_2) \\
&\approx s_2s_1ets_2etes_2ts_2s_1s_2s_1es_1s_2 &&\hspace{-3cm}\text{by (R11) and (R15)}\\
&\approx s_2s_1ets_2es_2ts_2s_2s_1s_2es_1s_2 &&\hspace{-3cm}\text{by (R13) and (R14)}\\
&\approx s_2s_1etes_2s_2ts_1s_2es_1s_2 &&\hspace{-3cm}\text{by (R11) and (R16)}\\
&\approx s_2s_1etes_2s_1s_2 &&\hspace{-3cm}\text{by (R11), (R14), (R15) and (R16)}\\
&\approx s_2s_1es_2s_1s_2 &&\hspace{-3cm}\text{by (R14),}
\end{align*}
with a similar calculation giving $\ep_3 \tau_{32} \ep_2 \tau_{21} \ep_1 \tau_{13} \ep_3 \approx s_2s_1es_2s_1s_2$. \epf

\ms
\begin{lemma}\label{lem:esse}
If $1\leq r\leq n$ and $1\leq k\leq n-1$, then $\ep_rs_k$ and $s_k\ep_r$ are both $\approx$-equivalent to elements of $\im\psi$.
\end{lemma}

\pf Since $s_k\ep_r\approx s_k\ep_rs_ks_k\approx \ep_{r\sb_k}s_k$, by (R11) and Lemma \ref{lem:ses}, it suffices to show that
\[
\ep_rs_k \approx 
\begin{cases}
\ep_r\tau_{r,r-1}\ep_{r-1}&\text{if $k=r-1$}\\
\ep_r\tau_{r,r+1}\ep_{r+1} &\text{if $k=r$}\\
\ep_r\tau_{rk}\ep_k\tau_{k,k+1}\ep_{k+1}\tau_{k+1,r}\ep_r &\text{otherwise.}
\end{cases}
\]
We just treat the case in which $k\not\in\{r-1,r\}$, the others being easier.  As above, we may assume that $k=1$ and $r=3$.
%By conjugating by a word $w\in S^*$ with $1\wb=k$ and $2\wb=k+1$, we may again assume that $k=1$.  But then 
%$\ep_2s_1=(s_1es_1)s_1\approx s_1e$ and 
%$\ep_2\tau_{21}\ep_1 = (s_1es_1)(s_1ts_1)e \approx s_1ete\approx s_1e\approx s_1es_1s_1=\ep_2s_1$, by (R11), (R14) and (R15), giving the $k=r-1$ case.  The $k=r$ case is similar.  For the third case, we may assume that $r=3$, and 
%By (R13) and the calculation from the proof of Lemma \ref{lem:reduction1}, we have $\ep_3 \tau_{31} \ep_1 \tau_{12} \ep_2 \tau_{23} \ep_3 \approx s_2s_1es_2s_1s_2 \approx s_2s_1es_1s_2s_1 = \ep_3s_1$. \epf
Here we have $\ep_3 \tau_{31} \ep_1 \tau_{12} \ep_2 \tau_{23} \ep_3 \approx s_2s_1es_2s_1s_2 \approx s_2s_1es_1s_2s_1 = \ep_3s_1$, by (R13) and the calculation from the proof of Lemma \ref{lem:reduction1}. \epf
%, using (R13) and the calculation from the proof of Lemma \ref{lem:reduction1}. \epf

\ms
\begin{lemma}\label{lem:tsst}
If $\oijn$ and $1\leq k\leq n-1$, then $\tau_{ij}s_k$ and $s_k\tau_{ij}$ are both $\approx$-equivalent to elements of $\im\psi$.
\end{lemma}

\pf Again, it suffices to do this for just $\tau_{ij}s_k$.  By (R11) and Lemmas \ref{lem:sts} and \ref{lem:reduction1}, we have $\tau_{ij}s_k \approx \tau_{ij}\ep_i\tau_{ij}s_k \approx \tau_{ij}\ep_is_ks_k\tau_{ij}s_k \approx \tau_{ij}(\ep_is_k)\tau_{i\sb_k,j\sb_k}$, and the result now follows from Lemma \ref{lem:esse}. \epf

%\pf Again, it suffices to do this for just $\tau_{ij}s_k$.  Since $\tau_{ij}s_k \approx \tau_{ij}\ep_i\tau_{ij}s_k \approx \tau_{ij}\ep_is_ks_k\tau_{ij}s_k \approx \tau_{ij}(\ep_is_k)\tau_{i\sb_k,j\sb_k}$, by (R11) and Lemmas \ref{lem:sts} and \ref{lem:reduction1}, the result now follows from Lemma \ref{lem:esse}. \epf

\ms
\begin{lemma}\label{lem:reduction2}
If $w\in(S\cup\{e,t\})^*\sm S^*$, then $w$ is $\approx$-equivalent to an element of $\im\psi$.
\end{lemma}

\pf Put $\Si=(E\cup T)\psi=\set{\ep_r}{1\leq r\leq n}\cup\set{\tau_{ij}}{\oijn}$, noting that $\im\psi=\la\Si\ra$.  Since $e=\ep_1\in\Si$ and $t\approx s_1ts_1=\tau_{12}\in\Si$, it suffices to show that every element of $\la\Si\cup S\ra\sm S^*$ is $\approx$-equivalent to an element of $\la\Si\ra$.  With this in mind, let $w\in\la\Si\cup S\ra\sm S^*$, and write $w=x_1\cdots x_k$, where $x_1,\ldots,x_k\in\Si\cup S$.  Denote by $l$ the number of factors $x_i$ that belong to~$S$.  We proceed by induction on $l$.  
If $l=0$, then we already have $w\in\la\Si\ra$, so suppose $l\geq1$.  Since $w\not\in S^*$, there exists $1\leq i\leq k-1$ such that either (i) $x_i\in S$ and $x_{i+1}\in\Si$, or (ii) $x_i\in\Si$ and $x_{i+1}\in S$.  In either case, Lemmas \ref{lem:esse} and \ref{lem:tsst} tell us that $x_ix_{i+1}\approx u$ for some $u\in\im\psi=\la\Si\ra$.  But then $w\approx (x_1\cdots x_{i-1})u(x_{i+2}\cdots x_k)$, and we are done, after applying an induction hypothesis (noting that $(x_1\cdots x_{i-1})u(x_{i+2}\cdots x_k)$ has $l-1$ factors from $S$). \epf

We may now prove the second main result.

\pfthm{thm:Pn} It remains to show that $\ker\Phi\sub{\approx}$, so suppose $w_1,w_2\in(S\cup \{e,t\})^*$ are such that $\wb_1=\wb_2$.  If $\wb_1\in\S_n$, then $w_1,w_2\in S^*$, so $w_1\approx w_2$, using only relations (R11--R13).  So suppose $\wb_1\not\in\S_n$.  It follows that $w_1,w_2\in(S\cup \{e,t\})^*\sm S^*$.  So, by Lemma \ref{lem:reduction2}, $w_1\approx u_1\psi$ and $w_2\approx u_2\psi$ for some $u_1,u_2\in(E\cup T)^+$.  We then have
$\ub_1=\overline{u_1\psi}=\wb_1=\wb_2=\overline{u_2\psi}=\ub_2$,
so that $u_1\sim u_2$, by Theorem \ref{thm:PnSn}.  Lemma \ref{lem:reduction1} then gives $u_1\psi\approx u_2\psi$, so that $w_1\approx w_2$. \epf

%\ms
%\begin{rem}
%Several other presentations for $\P_n$ were given in \cite{JEgrpm}, as well as presentations for the corresponding partition \emph{algebras}.
%\end{rem}

\footnotesize
\def\bibspacing{-1.1pt}
\bibliography{biblio}
\bibliographystyle{plain}
\end{document}